\documentclass[oneside,reqno,12pt]{amsart} 
\usepackage{amsmath}
\usepackage{amssymb}
\usepackage{amsthm}\usepackage{bbm}
\usepackage{amsxtra}
\usepackage{t1enc}
\usepackage{color}
\usepackage{epsf}
\usepackage{amssymb}
\usepackage{verbatim}
\newcommand\NoBlackBoxes{\global\overfullrule0pt}
\NoBlackBoxes
\theoremstyle{plain} 
\newtheorem{theorem}{Theorem} 

\newtheorem{lemma}[theorem]{Lemma}

\newtheorem{rem}[theorem]{Remark}

\newtheorem{corollary}[theorem]{Corollary}
\def\4{\kern1pt}

\def\6{\vphantom0}

\def\8{\kern-10pt}
\def\7#1{_{(#1)}}

\def\ov{\overline}

\theoremstyle{definition}
\newtheorem*{assumption (H1)}{Assumption (H1)}
\newtheorem*{assumption (H2)}{Assumption (H2)}

\theoremstyle{remark}

\numberwithin{equation}{section}
\numberwithin{theorem}{section}
\makeatletter
\let\serieslogo@\relax
\let\@setcopyright\relax

\def\speciallabelmark#1{\def\@currentlabel{#1}}
\makeatother

\newcommand{\epspot}{-K}

\newcommand{\abs}[1]{\lvert#1\rvert}
\newcommand{\norm}[1]{\lVert#1\rVert}
\newcommand{\eps}{{\varepsilon}}

\newcommand{\Cal}{\mathcal}

\newcommand{\Z}{{\mathbb Z}}
\newcommand{\Zd}{{\mathbb Z}^d}
\newcommand{\R}{{\mathbb R}}
\newcommand{\Rd}{{\mathbb R}^d}
\newcommand{\N}{{\mathbb N}}

\newcommand{\C}{{\mathbb C}}
\newcommand{\Q}{{ Q}}

\newcommand{\I}{{\mathbb I}}

\newcommand{\defi}{\;\overset {\text{\rm def}} = \;}
\DeclareMathOperator{\volu}{\operatorname{vol}}
\DeclareMathOperator{\Rea}{\operatorname{Re}}
\DeclareMathOperator{\Ima}{\operatorname{Im}}
\newcommand{\beqa}{\begin{eqnarray}}
\newcommand{\beqan}{\begin{eqnarray*}}
\newcommand{\eeqa}{\end{eqnarray}}
\newcommand{\eeqan}{\end{eqnarray*}}

\newcommand{\bw}{\noindent \textbf{Proof}.$\quad$}
\newcommand{\bwthm}[1]{\noindent \textbf{Proof of #1}.$\quad$}
\newcommand{\bwend}{{\hfill $ \square $}\\[2mm]}
\newcommand{\nn}{\nonumber}
\newcommand{\ffrac}[2]{\raise.5pt\hbox{\small$\4\displaystyle\frac{\,#1\,}{\,#2\,}\4$}}
\newcommand{\ovln}[1]{\,{\overline{\!#1}}}
\newcommand{\imbound}{\frac{1}{r}}
\newcommand{\imboundtwo}{\frac{2}{\pi r}}
\newcommand{\boundtappot}{-1}
\newcommand{\boundtap}{r^{\boundtappot}}
\newcommand{\volz}[2]{V^{\Z}_{#1 , #2}(r;a,b,M)}
\newcommand{\volr}[2]{V^{\R}_{#1 , #2}(r;a,b, M)}
\newcommand{\vzapprox}[4]{V^{\Z}_{#1 , #2}(r;#3,#4,M)}
\newcommand{\vrapprox}[4]{V^{\R}_{#1 , #2}(r;#3,#4,M)}
\NoBlackBoxes

\begin{document}

\title[Values of Special Indefinite Quadratic Forms]{Values of Special Indefinite Quadratic Forms${}^*$}

\author[G. ~Elsner]{Guido Elsner}
\address{\scriptsize Fac. of Mathematics, 
Univ. Bielefeld, P.O.Box 100131, 33501 Bielefeld, Germany}

\email{GUIDE@math.uni-bielefeld.de}

\thanks{${}^*$ Research supported by the DFG, CRC 701}


\date{February, 2007}

\begin{abstract}
For special $d$-dimensional hyperbolic shells $E$ with $ d\geq 5$ we
show that the number of lattice points in $E$ intersected with a
$d$-dimensional cube $C_r$ of edge length $r$, can be approximated
by the volume of $E\cap C_r$, as $r$ tends to infinity, up to an
error of order ${\mathcal O}(r^{d-2})$. We generalize results and
techniques, used by F. G\"otze (2004), to a large class of {\em
indefinite} quadratic forms and we provide explicit bounds for the
errors in terms of certain Minkowski minima related to these
quadratic forms. Furthermore, we obtain, as in the positive definite
case, a result for multivariate diophantine approximation and for the maximal gap between values of such indefinite forms.
\end{abstract}

\maketitle

\section{Introduction}
 
Let $\Q $ denote a $d$-dimensional quadratic form. For $a,b \in \R$ we consider the set $ E$ of points in the $d$-dimensional Euclidean space, for which $Q$ takes values between $a $ and $b$. In case that the quadratic form $ \Q[\4x\4] $ is positive definite, $E$ is an elliptic shell, but in this paper we will investigate indefinite forms and hence, $E$ is a hyperbolic shell. \\
For a (measurable) set $B\subset \Rd$ the lattice volume of $B$ is
the number of lattice points in $B$ (formally $\volu_{\mathbb{Z}} B
\defi \#\bigl(B\cap \Zd\bigr)$) and $\volu B$ denotes the Lebesgue
measure of $B$. For the hyperbolic shell $E$ we want to approximate its
lattice point volume by the Lebesgue volume. We want to investigate
this approximation by estimating a {\it relative} lattice point rest
of {\em large} parts of the hyperbolic shell $E$. Therefore we consider
for $r>0 $ the $d$-dimensional  cube $C_r $ with edge length $r$ and
intersect the cube $C_r$ with the hyperbolic shell $E$. The
 {\it relative}
lattice point rest of $E\cap C_r$ is now defined by \beqa
\Delta \defi
\Bigl|\ffrac{\volu_{\mathbb{Z}} (E\cap C_r)-\volu (E\cap C_r)}{\volu (E\cap C_r)}\Bigr|.\nn
\eeqa

We will show for special indefinite forms $\Q$, that $\Delta =
\Cal{O}(1)$ as $r\to \infty$ (Theorem \ref{1.1})
and that even  $\Delta = o(1)$  as $r\to
\infty $ provided that $\Q $ is irrational (Theorem \ref{1.2}).
Recall that a quadratic form \,$\Q[\4x\4]$ \,and the corresponding
operator $\Q$ with non-zero matrix \,$\Q=(q_{ij})$, \,$ 1 \leq i, j
\leq d$, \,is called {\it rational} if there exists a real number
\,$\lambda \neq 0$ \,such that the matrix \,$\lambda\4 \Q$ \,has
integer
 entries only; otherwise it is
called {\it irrational}.\\
Similar results for forms $\Q$ of  signature $(p,q)$  satisfying  ~$ \max (p,q) \ge 3$
have been proved
by Eskin, Margulis and Mozes in \cite{eskin-margulis-mozes:1998}.
These are quantitative versions of the well-known Oppenheim problem
concerning  the distribution  of values of ~$\Q[m]$, \,$ m \in \Zd$.
In 1929, Oppenheim (\cite{oppenheim:1929}, \cite{oppenheim:1931}) conjectured that if $d\geq 5 $ for an irrational non-degenerate quadratic form $Q$ the quantity $m(Q)\defi \inf\{\bigl|Q[m]\bigr| : m \in \Zd, m\neq 0\}$ equals zero. In the rational case this was known by Meyer's Theorem (see \cite{cassels:1978}). Later it was conjectured that even for $d\geq 3$ and $Q$ irrational the equality $m(Q) = 0 $ holds (for irrational diagonal forms this was suspected in \cite{davenport-heilbronn:1946} and it is not true in dimensions 3 and 4 without the assumption of irrationality). The different approaches to this and related problems involve various mathematical methods from analytic number theory, from ergodic theory, from representation theory of Lie groups, reduction theory and from the geometry of numbers. In \cite{margulis:1987} Margulis  established the Oppenheim conjecture in dimensions $d \geq 3$, as stated by
Davenport and Heilbronn for $d\geq 5$. In his seminal work he proved that the set of values of $Q$ at lattice points is dense in $\R$. Quantitative versions of this problem were later on developed by Dani and Margulis (\cite{dani-margulis:1993}) and Eskin, Margulis and Moses (\cite{eskin-margulis-mozes:1998}). They consist of {\em quantative} bounds on the ratio between the lattice point volume and the Lebesgue volume of the set of points in the cube $C_r$, where the quadratic form takes values in a small interval. 
The quantitative bounds provided in these results yield the asymptotic number of points in these regions as a polynomial in $r$ up to a non-effective error term tending to zero in proportion to the leading term.
The estimates thus obtained are implicit, since they do not provide explicit bounds in terms of diophantine approximations of irrational coefficients of the form.
For a detailed
discussion of results on these problems
 by Oppenheim, Heilbronn and Davenport and others, see \cite{margulis:1997}. In \cite{bentkus-goetze:1999} Bentkus and G\"otze proved explicit error bounds in the quantitative Oppenheim problem for the elliptic shell
as well as for hyperbolic shells for $d\ge 9$  by a common approach. They provide more explicit bounds (in terms of diophantine approximation) for {\em distribution functions} of the values of the quadratic form on $C_r$, whereas the direct application of the previous methods seems to be restricted to the case of the concentration in compact intervals. \\
In \cite{goetze:2004} G\"otze showed  that in the positive definite case for $d\geq 5$ the lattice point rest is of order $\Cal{O}(r^{d-2})$ for arbitrary forms, and of order $o(r^{d-2})$ if the form is irrational. These results refine earlier bounds of the same order for dimensions $d\geq 9$ (see also \cite{goetze:2004} for the history of such estimates and further references).\\[1mm]
In the present paper we apply techniques of
\cite{goetze:2004} to special {\em indefinite} forms and we obtain
explicit bounds in terms of certain Minkowski minima of convex
bodies related to these quadratic forms.
Adapting these techniques, the main problem consists of the estimation of the difference between the lattice point and the Lebesgue volume by an integral of generalized theta functions. In order to achieve such an estimate, we develop tools, different from those in \cite{goetze:2004}, which involve adjustable smooth approximations of the indicatior functions of the hyperboloid and of the cube $C_r$. The bound given by an integral of theta functions does {\em not} use the special structure of the indefinite forms under consideration. Furthermore, a careful modification of the arguments in \cite{goetze:2004} even leads to a bound in terms of the Minkowski minima mentioned above, which holds for {\em any} indefinite form. The special structure of the forms is only used when we estimate the appearing functions of Minkowski minima by adapting the techniques of \cite{goetze:2004} to the indefinite case.
As in the positive definite case we show that in the irrational case the maximal gap between successive values of the quadratic form at lattice points converges to 0 as $ r $ tends to infinity (Corollary \ref{corgaps}). 
Furthermore, we extend the results of Bentkus and G\"otze (\cite{bentkus-goetze:1999}) on distribution functions for values of  quadratic forms to dimensions including 5 up to 8 (Theorem \ref{thdistr}).
In addition, we obtain a result for multivariate diophantine approximations for these special indefinite forms (Theorem \ref{diophant}).\\[1mm] This paper is organized as
follows: In the second section, we state the two main results about
the asymptotics of the relative lattice point rest and derive two important
corollaries concerning gaps between values of the quadratic form and
concerning multivariate diophantine approximations. Furthermore, we
give explicit quantitative bounds for the relative lattice point
rest. In the third section, we prove the results of the second section. In the
fourth section, we collect auxiliary results (e.g. from geometry of
numbers, metric number theory, theory of theta functions), which are
used in the proofs of the theorems.\\[2mm] 
{\em Acknowledgements}\\
I would like to thank Prof. Dr. Friedrich G\"otze for drawing my attention to this topic, for various fruitful discussions and many valuable suggestions. Furthermore, I am grateful to the DFG-CRC 701 for financial support. This paper is a part of my PhD thesis \cite{elsner:2006}.\newpage


\section{Results}

Let \,$\mathbb{R}^d$, \,$ 1 \le d < \infty$, \,denote the $d$-dimensional
Euclidean space with scalar product $\langle\cdot, \cdot\rangle$
and norm  $\lvert\cdot\rvert$ defined by
\,$\abs{ x }^2 =
\langle x,x \rangle = x_1^2 + \dots+ x_d^2$
for $ x =(x_1, \dots, x_d) \in
\mathbb{R}^d$. \,Let $\mathbb{Z}^d$ denote the standard lattice of points
with integer coordinates in $\mathbb{R}^d$.\\
Consider the quadratic form
$$
\Q[\4x\4] \defi \langle\4 \Q\4 x ,x\4\rangle ,
 \quad\text{for} \quad x \in \Rd,
$$
where $\Q:\Rd\to \Rd $ denotes a symmetric linear operator
in ${\rm{GL}} (d, \mathbb{R})$ with eigenvalues,
 say, ~${q_1,\dots, q_d}$. We write
\beqa
 q_0\defi\min\limits_{1\leq j\leq d}\, \abs{q_j},
\quad\quad\quad q\defi\max\limits_{1\leq j\leq d}\, \abs{q_j},\quad\quad\quad \bar{q}\defi\max\left\{q_0^{-1} ; q\right\}.  \label {eq:1.1}
\eeqa
In the sequel we always assume that the form is non-degenerate, that is, that $q_0>0$. \\[2mm]
We say that a quadratic form $\Q $ is of {\em block-type}, if and only if we can write $\Q = \Q^+ -\Q^-$, where $\Q^+ $ and $\Q^-$ are positive definite quadratic forms, $\Q^+[x]$ only depends  on the first $d_1$ coordinates of $ \Rd$ and $\Q^-[x]$  on the $d-d_1$ remaining ones only.\\[2mm]
We define for $ a, b \in \R,$ with $ a\leq b $ and for $M \in \Rd$ the sets
\beqa
E_{a,b;M} \defi \bigl\{ x \in \Rd \; : \; a\leq \Q[\4x - M\4] \leq b \bigr\} .
\eeqa
Note, that if the quadratic form $ \Q[\4x\4] $ is positive definite,
then $ E_{a,b;M}$  is an elliptic shell.

Recall that a quadratic form \,$\Q[\4x\4]$ \,and the
corresponding operator $\Q$
with non-zero matrix \,$\Q=(q_{ij})$, \,$ 1 \leq i, j \leq d$, \,is called
{\it rational} if there exists a real number \,$\lambda \neq 0$ \,such
that the matrix \,$\lambda\4 \Q$ \,has integer
 entries only; otherwise it is
called {\it irrational}.

For $r>0 $ we set $C_r \defi \{ x\in \Rd \; : \; \abs{x}_{\infty} \leq r\}$, where $\abs{\cdot}_{\infty} $ denotes the maximum norm on $\Rd$, and
\vskip-2mm \beqa  H_{r,M}\defi H_{r,M}^{a,b}\defi E_{a,b;M}\cap C_r.\label{defHrab}
\eeqa

For any (measurable) set $B \subset \mathbb{R}^d$ let $\volu B$
denote the Lebesgue measure of $B$ and $\volu_{\mathbb{Z}} B$ its
lattice volume, that is the number of lattice points in \,$B \cap
\mathbb{Z}^d$.
We want to investigate the approximation of the
lattice volume of $H_{r,M}$ by the Lebesgue volume. Therefore we  estimate the following {\it relative}
lattice point rest of {\em large} parts of hyperbolic shells $ H_{r, M}$, $M \in \Rd$, $r$ large.
\pagebreak\\

We define
\beqa
\Delta(r,M) \defi
\Bigl|\ffrac{\volu_{\mathbb{Z}} H_{r,M}-\volu H_{r, M}}{\volu H_{r,M}}\Bigr|.
\label{defDeltarM}
\eeqa
\vskip3mm

The two main results of this part of the paper are the following

\begin{theorem}{\label{1.1}} For a non-degenerate, $d$-dimensional, block-type form $\Q$, $d\geq 5$ and  all $M \in \mathbb{R}^d$  it holds
\beqa
\Delta(r,M) = {\Cal O}(1),
\qquad \text{ as } r\to\infty.\label{eq:gen}
\eeqa

\end{theorem}

The estimate of Theorem \ref{1.1} refines an explicit bound
of order $\Cal O(1)$
obtained for dimensions $d\ge 9$ in \cite{bentkus-goetze:1997} 
for {\it arbitrary} ellipsoids and  in \cite{bentkus-goetze:1999} for {\it arbitrary} hyperbolic shells. Since this bound is optimal in the case of positive definite forms (\cite{goetze:2004}, p. 196), the bound in Theorem \ref{1.1} is also optimal for block-type forms.\\[2mm]
In case that $\Q$ is {\it irrational} Theorem \ref{1.1} can be improved.

\begin{theorem}{\label{1.2}}
 For a non-degenerate $d$-dimensional block-type form $\Q$,\linebreak $d\geq 5$ and all $M \in \Rd$ it holds
\beqa \label{eq:irr}
\Delta(r,M ) =  o(1),
\quad \mbox{ provided that
 }\, \Q \, \mbox{ is irrational}.
\eeqa
\end{theorem}
\vskip2mm

For {\it irrational} forms and dimension
 $d \ge 9$ the bound of Theorem~\ref{1.2} has been already proved in
\cite{bentkus-goetze:1999}. We should remark again, that the bounds of both theorems are explicit and effective.
\begin{rem}\label{remark}
For $M \in \mathbb{Q}^d$ the condition $\Delta(r,M ) =  o(1)$ implies that $\Q$ is irrational.
\end{rem}
Using Theorem \ref{1.2} we can derive easily a Corollary about gaps between values of block-type forms:\\[3mm]
For a positive definite quadratic form, Davenport and  Lewis
\cite{davenport-lewis:1972} conjectured in 1972, that the distance
between successive values $v_n$ of the quadratic form ~$\Q[x]$ on
~$\Zd$ converges to zero as $n\to \infty$, provided that
the dimension $d$ is at least five and $\Q$ is irrational. This conjecture was proved by G\"otze \cite{goetze:2004}. Now we can derive an analog result for irrational block-type forms and dimension $d\geq 5$.\\[3mm]
For a vector $M\in \Rd$ let
\beqa
V(r)\defi \bigl\{ \Q[x-M]\; : \;x \in \Zd \cap C_r\}
\eeqa denote the set of values
of $\Q[\4x-M\4]$,  for lattice points $ x \in \Zd$ in a box of edge length $r$.\\ We define the maximal gap
between successive values as
\beqa
d(r) \defi \sup\limits_{u \in V(r)}\inf\bigl\{ v-u\; : \; v>u , v \in V(r)\bigr\}.\eeqa

\begin{corollary}{\label{corgaps}} For a non-degenerate $d$-dimensional block-type form $\Q$, 
\linebreak $d \ge 5$, the follwing holds:
\begin{itemize}
\item[(1)]
$\lim\limits_{r \rightarrow \infty}d(r) = 0,$ provided that $\Q$ is irrational.
\item[(2)] For $M\in \mathbb{Q}^d $ and $\Q $ rational we get  $
\lim\limits_{r \rightarrow \infty}d(r) > 0.$
\end{itemize}
\end{corollary}

The Theorems \ref{1.1} and \ref{1.2} follow from Theorem \ref{maindelta} below. Furthermore, in  Theorem \ref{maindelta} (combined with \eqref{vollowerbound} in the proof of the Theorems \ref{1.1} and \ref{1.2}),
 estimates of the remainder terms in \eqref{eq:gen}
 and \eqref{eq:irr} in terms of certain diophantine
properties  of $\Q$ will be given.\\[2mm]
In order to describe the explicit
bounds we need to introduce some more notations.
Let $\lvert (x,y) \rvert_{\infty}$ denote the maximum norm of a
 vector ~$(x,y)$ in $\Rd \times \Rd$. For any ~$t >0$ and ~$r\ge 2$ consider the
norm $F$ on $\Rd \times \Rd$ given by
\beqa
F(x,y) \defi \bigl\lvert\bigl( r\4(x + t\4\Q \4 y)\,,\, y\4r^{-1}\bigr) \bigr\rvert_{\infty}.
\eeqa
We introduce the so called Minkowski minima of the convex body $\{ F \le 1 \}$ as
\beqa\label{eq:norm0}
M_{1,t} = \inf\bigl \{F(m,n)\,:\,
(m,n) \in (\Zd\times \Zd) \setminus 0 \bigr\}
\eeqa
and we define in general
$M_{k,t}$ as the infimum of $\lambda > 0$
such that the set of lattice points  with norm less than ~$\lambda,$ that is $$\bigl\{(m,n) \in \Zd \times \Zd \,:\,
F(m,n) < \lambda\bigr \},$$ contains $k$ linearly independent vectors. By definition we have
~$r \4 M_{k,t} \ge 1$. For $d>4$ and $r\ge2$ we introduce
\beqa
\Gamma_{T,r} &\defi&
\inf\hskip5pt \bigl\{  r^d\4 M_{1,t}\cdots M_{d,t} \, : \, T^{-1/(d-4)}
\le  |t|\le T \bigr\},\label{eq:Gamma}\\[0.2cm]
\rho(r, \Q, T) &\defi&
\bar{q}^{d+1}T^{-\frac{1}{2}}+\bar{q}^{\frac{3d}{2}}\max\bigl\{\text{\small$\frac{2}{\pi r},\frac{\pi}{2q_0qr}$},T^{-\frac{1}{d-4}}\bigr\}  \nn\\ &&\quad \quad\quad\quad\quad+\quad
 \bar{q}^{d+2}\4  \Gamma_{T,r}^{-\frac{1}{2} +\frac{2}{d}}\4 \log\bigl(\bar{q}\4 T^{\frac{1}{2}}\4 \Gamma_{T,r}+1\bigr),
\label{defrhoT} \\
\rho(r, \Q) &\defi& \inf\limits_{T \geq 1} \Bigl\{r^{2-d} +q_0^{\frac{d}{2}}r^{2-\frac{d}{2}}+\bar{q}r^{2-\frac{d}{2}}(1+\log r) +\rho (r,\Q, T)  \Bigr\}\nn\\[-1mm] &&.\label{defrho} \eeqa For any
fixed $T>1$ and irrational $\Q$ it is shown in Lemma
\ref{irrational}\/ that \beqa \lim\limits_{r\to \infty} \Gamma_{T,r}
= \infty, \eeqa with a speed depending on the diophantine properties
of \,$\Q$. \,This implies that \beqa\label{eq:conv} {\lim\limits_{r
\to \infty} \rho(r, \Q) =0}. \eeqa
 With these notations we may state a Theorem providing quantitative bounds for the difference between the volume and the lattice point volume of a hyperbolic shell.
\begin{theorem}\label{maindelta}Let $\Q$ denote a non-degenerate $d$-dimensional block-type form, $d\geq 5$, and $ M \in \Rd$. Furthermore, let $c(Q,M)>0$ be defined as in Theorem \ref{main} below and $K=K(d)$ is chosen according to \eqref{taperingfunction}. Then there exist
 constants \,$c_j>0$, \,$j=1,2$, \,depending on $d$ only
and a constant $r_0 = r(\Q, M, a, b)>0$ such that, for any $r \ge r_0,$\vskip3mm
\begin{itemize}
\item[(1)]$ \bigl|\volu_{\mathbb{Z}}H_{r,M}-\volu H_{r,M}\bigr|$
$$\quad\quad \le \quad c_1\cdot r^{d-2}\cdot \Bigl((b-a+1)\bar{q}^{d}q^{-1} + c(Q,M)\bar{q}^{d+1}(\log q +1)+1\Bigr).$$
\item[(2)]$
\bigl|\volu_{\mathbb{Z}}H_{r,M}-\volu H_{r,M}\bigr| \;\;\le
\;\; c_2\cdot r^{d-2} \cdot \Bigl(
(b-a)\,\bar{q}^{d}q^{-1}r^{-\frac{1}{K}} $  $$\quad\quad\quad\quad + \quad
(b-a)\,\bar{q}^{d+1}q^{-1}\bigl(\abs{M}+2q^{-\frac{1}{2}}\frac{\abs{a}+\abs{b}}{r}\bigr)r^{-1} \, +
\, c(Q,M)\cdot\rho(r,\Q )\Bigr),$$ where $ \lim\limits_{r\to \infty}
\rho(r, \Q) =0$,\;\; provided that\;$  \Q $ is irrational.
\end{itemize}
\end{theorem}

Note that the  summand $\rho(r,\Q)\4 r^{d-2}$ in the bound in Theorem \ref{maindelta}
 is at least of order
${\Cal O} (r^{d/2}\4 \log r).$ It may be indeed  of  this order
 since  $r\4 M_{j,t} \ll_d r$ shows that the  maximal
value of $\Gamma_{T,r}$ is of order $\Cal O (r^{d})$ and we may choose
$T= {\Cal O}(r^{\bar{\beta}})$ with $\bar{\beta} > 0$ sufficiently large.\\
Note that  an  error bound  of order $r^{d/2+\eps}$ has been
proved by Jarnik \cite{jarnik:1928}
for {\it diagonal} $\Q$ = diag$(s_1,..,s_d),\, s_j > 0$  for Lebesgue {\it almost all} coefficients
$s_j.$\\[2mm]
The proof of Theorem \ref{maindelta} is based, roughly speaking, on
an 'continuous' approximation of $
\bigl|\volu_{\mathbb{Z}}H_{r,M}-\volu H_{r,M}\bigr|$ by an integral
over generalized theta functions. We will derive bounds for parts of
this integral, which use the distribution of the first Minkowski
minimum $M_{1,t}$. We investigate this distribution using results
from metric number theory. As a consequence of this investigation,
we also get a result for multivariate diophantine
approximation:\\[2mm]
For a vector ~$x\in \Rd$ let ~$\| x \| \defi \inf \limits_{m \in
\Zd} \vert x-m \rvert_{\infty} $  denote the error of an integer
approximation. For
 real numbers  $t>0$, $\nu > 1$ we introduce
\beqa D(t,\nu)= \nu \4\min \, \bigl\{ \| t\4 \Q\4 n\|  \,:\, n \in
\Zd, \, 0 < \lvert n\rvert_{\infty} \le \nu  \bigr\},
\label{defD(t,v)}\eeqa and let $\lambda$ denote the Lebesgue
measure. Then we have
\begin{theorem} \label{diophant}
Assume that $\Q$ is a symmetric, non-degenerated block-type form,
which is normalized such that $q_0=1.$ Then there exists a constant
$c(d)>1$ depending on $d$ only such that for any $r\ge 1$ and any
interval $[\kappa,\xi]$ satisfying $0<\xi-\kappa < 1$ the following
inequalities hold \beqa\label{eq:metric1}
\lambda\{t\in [\kappa,\xi] \, :\, M_{1,t} \le \tau\}& \le& c(d)\bigl( q \4 \tau^2 \4(\xi-\kappa) + \tau\4 r^{-1}\bigr),\\[0.1cm]
\label{eq:metric2} \sup\limits_{t\in [\kappa,\xi]} \,M_{1,t}&\ge&
\min\{\tau_Q, r\4(\xi-\kappa)\},
\\[0.1cm]\label{eq:metric3}
\sup\limits_{t\in [\kappa,\xi]} D(t,\nu)  &\ge&  \min\{ \tau_Q,
\nu\4 (\xi-\kappa)/2 \}, \eeqa for any $ \nu \ge \tau_Q$, where $
\tau_Q\defi \left(\frac{c(d)+2}{2c(d)q}\right)^{1/2}$.
\end{theorem}

Refining the proofs, we may extend Theorem \ref{1.1} and \ref{1.2} to include the case $a = - \infty$, i.e. the case of distribution functions. This partially extends a result obtained by Bentkus and G\"otze in \cite{bentkus-goetze:1999} to the dimensions including 5 up to 8. 

\begin{theorem}\label{thdistr}
For a non-degenerate, $d$-dimensional, block-type form $\Q$, $d\geq 5$ and  all $M \in \mathbb{R}^d$ we set
$$F_{r,M}(b)\defi \big\{x \in \Rd \; : \; \Q[\4x - M\4] \leq b, \;\;\abs{x}_{\infty}\leq r\big\}.$$

 Then for the corresponding relative lattice point remainder holds
$$
\Bigl|\ffrac{\volu_{\mathbb{Z}} F_{r,M}(b)-\volu F_{r, M}(b)}{\volu F_{r,M}(b)}\Bigr| = \left\{\begin{array}{ll} o(1),& \text{ provided } \Q \text{ is irrational,}\\[3mm] \Cal{O}(1),& \text{ otherwise},\end{array}\right.$$ as $ r \rightarrow \infty.$
\end{theorem}

\section{Proofs}
First we deduce Theorem \ref{1.1} and \ref{1.2} from Theorem \ref{maindelta}:\\

\bwthm{Theorem \ref{1.1} and \ref{1.2}}\\
By Lemma \ref{vol3} we obtain for $M = (M_1,...,M_d)$ and for $r$ large
\beqa
\volu \, H_{r, M} &\gg_d & (b-a)q^{-d/2}\left(q_0^{d/2} + r^{-1}\bigl| \bigl(\abs{q_1}^{\frac{1}{2}}M_1 ,...,\abs{q_d}^{\frac{1}{2}}M_d\bigr)\bigr|\right)^{d-2}r^{d-2}\nn\\ &\gg_d& (b-a)q^{-d/2}q_0^{d(d-2)/2}r^{d-2}.\label{vollowerbound}
\eeqa
Dividing the inequalities in Theorem \ref{maindelta} (1) in the general case (resp. Theorem \ref{maindelta} (2) in the irrational case) by $ \volu \, H_{r,M}$, the estimate \eqref{vollowerbound} completes the proof of Theorem \ref{1.1} (resp. of Theorem \ref{1.2}).
 \bwend
\bwthm{Corollary \ref{corgaps}}\\
If $\Q$ is irrational, Theorem \ref{1.2} implies, that for any $a,b
\in \R$ \beqa \left|\frac{\volu_{\Z} H^{a,b}_{r,M}}{\volu_{\R}
H^{a,b}_{r,M}}-1\right| \rightarrow 0 \;\;\text{ as }\;\;
r\rightarrow \infty. \eeqa Hence,   $H_{r,M}^{a,b}\cap \Z^d \neq
\emptyset$ for all $a,b \in \R$ if $r$ is sufficiently large. This
implies that $\lim\limits_{r \rightarrow \infty} d(r) =0$ and the proof of part (1) is completed. \\
If  $\Q$ is rational, there exists a real number $\lambda>0$, such that $\lambda\Q$ has integer entries only. For $M\in \mathbb{Q}^d$ there exists a $\mu \in \Z$, $\mu \neq 0,$ such that $\mu M\in\Zd$. Hence, it holds that $\Q[m-M]\in \lambda^{-1}\mu^{-2}\Zd$ for all $m\in \Zd$. Therefore $d(r)\geq \lambda^{-1}\mu^{-2}>0$ for all $r\geq 1$, which proves part (2).
\bwend
We should remark, that  by using \eqref{2.3main} and \eqref{vollowerbound} one can obtain explicit bounds for $d(r)$ in terms of $r$ and $\rho(r,\Q)$, representing diophantine properties of $\Q$.\\[5mm]
\bwthm{Remark \ref{remark}}
Analyzing the proof of Corollary \ref{corgaps} (1) we recognize that $\Delta(r,M) = o(1)$ already implies $\lim_{r\to\infty} d(r) = 0$. Under the assumption $M\in\mathbb{Q}^d$ the condition that $Q$ is rational yields by Corollary \ref{corgaps} (2) that $\lim_{r\to\infty} d(r) > 0$. Thus, for $M \in \mathbb{Q}^d$ the irrationality of $\Q$ follows from $\Delta(r,M) = o(1)$. 
\bwend

The first step in proving Theorem \ref{maindelta} is to analyze
smooth approximations of the
lattice volume of $H_r$:\\[2mm]
For $a,b \in \R$ and a smoothing parameter $w> 0$ we define $ g_{a,b,w}: \R \rightarrow [0,1]$ by \beqa
g_{a,b,w}(x)\defi  \frac{1}{w}\bigl( (b+w-x)_+ - (b-x)_+ - (a-x)_+ +
(a-w-x)_+\bigr).\label{defg} \eeqa This function $g_{a,b,w}$ is a
linear continuous approximation of the indicator function
$I_{[a,b]}$ of the interval $[a,b]$. By Lemma \ref{pospartintegral} we may rewrite $g_{a,b,w} $ as follows
\beqa
g_{a,b,w} (x)&= &\frac{1}{2\pi i} \int\limits^{\beta+ i\, \infty}_{\beta-i\, \infty}\!
e^{(b+w-x)z} - e^{(b-x)z} - e^{(a-x)z} + e^{(a-w-x)z} \,\ffrac {dz} {wz^2} \nonumber
\\ &=& \frac{1}{2\pi i} \int\limits^{\beta+ i\, \infty}_{\beta-i\, \infty}\exp\left[-xz\right]\cdot h_{a,b,w}(z)\frac{dz}{z},\label{g-integral}
\eeqa
where $ h_{a,b,w}(z)\defi \ffrac{\exp[wz]-1}{wz}\cdot\Bigl(\exp[bz] - \exp [(a-w)z] \Bigr)$.\\

Using $g_{a,b,w} $ we construct a continuous approximation $\volz{w}{\eps}$ of the (monotone) lattice point counting function $ r \mapsto  \volu_{\mathbb{Z}}( H_{r,M})$ depending on two smoothing parameter $w>0$ and $\eps > 0$. Setting $Q_+\defi \big(Q^TQ)^{\frac{1}{2}}$, we define
\beqa
\volz{w}{\eps}\!\!\!\! &\defi& \!\!\!\!\sum_{x\in \Z^d} \exp\Bigl[-\frac{2}{r^2}\,Q_+[x]\Bigr]\,g_{a,b,w}\Bigl( Q[x-M] \Bigr)\, \chi_{\eps}\left(\frac{x}{r}\right) \label{volzdef}\\\text{and\hskip20mm} &&\nn\\
\volr{w}{\eps}  \!\!\!\! &\defi& \!\!\!\! \int\limits_{\R^d} \exp\Bigl[-\frac{2}{r^2}\,Q_+[x]\Bigr]\, g_{a,b,w}\Bigl( Q[x - M]\Bigr) \,\chi_{\eps}\left(\frac{x}{r}\right)dx\phantom{......}\label{volrdef}, \eeqa
where $\chi_{\pm \eps}$ is a function with the following properties:
\begin{itemize}
\item[(1)] For $u\in \R^d$ it holds $$\chi_{\pm\eps}(u) =  \left\{ \begin{array}{ll} \exp\left[2\cdot Q_+[u]\right],& \text{ if } |u|_{\infty }\leq \min \{1; 1\pm\eps\},\\0,& \text{ if }|u|_{\infty }> \max\{1;1\pm \eps\} .
\end{array}\right.$$
\item[(2)] There exists a constant $ c_1(Q,M) > 0$ such that for \beqa\label{defbarchi}\overline{\chi}_{\pm\eps}(x) \defi \chi_{\pm\eps}(x)\cdot\exp\bigl[\left\langle x , 2r^{-1}QM\right\rangle\bigr] \eeqa the following estimates hold for an appropriate $K=K(d) \in \N$:
\beqa  
&\text{(a)}& \int_{\Rd}\;\;\left|\widehat{\overline{\chi}}_{\pm\eps}(v)\right|dv \quad\ll_d \quad c_1(Q,M)\cdot \eps^{\epspot} ,\nn\\[-2mm] &&\label{taperingfunction}\\[-2mm]
&\text{(b)}& \int_{\{ |v|_{\infty} > d^{-\frac{1}{2}}r\}}\left|\widehat{\overline{\chi}}_{\pm\eps}(v)\right|dv \; \ll_d \; c_1(Q)\cdot \eps^{\epspot}\boundtap \;\text{ for all }\; r \geq 1.\nn
\eeqa
\end{itemize}

The existence of such a function $\chi_{\pm\eps}$ follows by standard arguments from Fourier analysis (cf. \cite{elsner:2006}, p. 27, Lemma 2.4.5).
Note, that  the function \beqa\psi_{r,\pm\eps}(x)&\defi& 
\exp\Bigl[-\frac{2}{r^2}\,Q_+[x]\Bigr]\chi_{\pm\eps}(\frac{x}{r})\label{defpsi_r}\eeqa approximates the indicator function $I_{\{\abs{x}_{\infty} \leq r\}}$ and hence the equations \linebreak
 $\volz{0}{0} = \volu_{\mathbb{Z}}(H_{r,M})$ and $\volr{0}{0} = \volu_{\R}( H_{r,M})$ are suggestive.\\

\bwthm{Theorem \ref{maindelta}}\\
For $M \in \Rd,\, 0<\eps \leq \frac{1}{4}$, there exists a constant $c = c(d) >0$ by Lemma \ref{approx1} such that \\

$\bigl|\volu_{\Z} H_{r,M} - \volu H_{r,M}\bigr| \;\;\leq\;\; \max
\bigl\{ \Delta_{-\eps} \; ; \; \Delta_{\eps}\bigr\}\;\; +  c\cdot
(b-a)q_0^{-d/2}q^{(d-2)/2}$\beqa \quad\quad\quad\quad\times\left(\eps + q_0^{-1/2}q^{1/2}|M|r^{-1} +
2q_0^{-1/2}(|a|+|b|)r^{-2} \right)\, r^{d-2},\label{eq1:maindelta}
\eeqa where $\Delta_{\pm\eps}$ is defined by using \eqref{defpsi_r} as follows
\beqa
\Delta_{\pm\eps} &\defi &\Bigl|\int_{\Rd} I_{H_{r,M}}(x)\psi_{r,\pm\eps}(x)dx - \sum_{x\in\Zd} I_{H_{r,M}}(x)\psi_{r,\pm\eps}(x)\Bigr|.
\label{defDeltaeps}
\eeqa\\
Hence, $\Delta_{\pm\eps}$ can be estimated by Lemma \ref{approx3} by \beqa
\Delta_{\pm\eps} &\ll_d &\max\limits_{\pm}\sup\limits_{{a' \in
[a-w;a+w]}\atop{b' \in
[b-w;b+w]}}\!\Bigr|\vrapprox{w}{\pm\eps}{a'}{b'} -
\vzapprox{w}{\pm\eps}{a'}{b'}\Bigr|\hskip2cm\nn\\&&\hskip3cm+ \; 8w
q_0^{-\frac{d}{2}}q^{\frac{d-2}{2}}
\bigl(1+\eps+q^{\frac{1}{2}}\frac{\abs{M}}{r}\bigr)^{d-2}r^{d-2}.\label{eq2:maindelta}
\eeqa
Collecting the estimates \eqref{eq1:maindelta} and \eqref{eq2:maindelta} we obtain for $w> 0, 0<\eps \leq \frac{1}{4}$\\

$r^{2-d}\cdot\bigl|\volu_{\Z} H_{r,M} - \volu H_{r,M}\bigr|
\;\;\ll_d \;\;(b-a)\,\bar{q}^{d}q^{-1}\eps$\\[-2mm]
\beqa & + &
(b-a)\,\bar{q}^{d+1}q^{-1}\bigl(\frac{\abs{M}}{r}+2q^{-\frac{1}{2}}\frac{\abs{a}+\abs{b}}{r^2}\bigr)\;+\;
w\bar{q}^{d}q^{-1}\bigl(1 +\eps +
q^{\frac{1}{2}}\frac{\abs{M}}{r}\bigr)^{d-2}\nn\\ & +&
\;\;\max\limits_{\pm}\;\sup\limits_{a' ,b'}\;
\Bigr|\vrapprox{w}{\pm\eps}{a'}{b'} -
\vzapprox{w}{\pm\eps}{a'}{b'}\Bigr|\cdot r^{2-d}. \label{eq3:maindelta} \eeqa

Choosing now $w = 1, \eps = \frac{1}{4} $ and $ r>r_0$ large
enough, \eqref{eq3:maindelta} and the result of the following
crucial Theorem \ref{main} (1) below yields (note that $d\geq 5$) \\

$ \bigl|\volu_{\Z} H_{r,M} - \volu H_{r,M}\bigr| $ \beqa &\ll_d&
(b-a+1)\bar{q}^dq^{-1}r^{d-2} + c(Q,M)\Bigl(1 + q_0^{-\frac{d}{2}}r^{\frac{d}{2}} + \bar{q}^{d+1}(\log q +1)\, r^{d-2}\Bigr)\nn\\
\nn 
&\ll_d & \Bigl((b-a+1)\bar{q}^{d}q^{-1} + c(Q,M)\bar{q}^{d+1}(\log q +1)+1\Bigr)r^{d-2}, \nn \eeqa for $r$ large enough.
This proves Theorem \ref{maindelta} (1).\\[3mm]
For proving Theorem \ref{maindelta} (2), we choose for an arbitrary $\alpha \in (0,1)$\beqa w =
T^{-1/2}, \, T\geq 1\;\text{  and } \;  \eps = \eps(r) \defi r^{-\alpha K^{-1}}
,\nn\eeqa where $ K = K(d)$ is chosen according to \eqref{defbarchi}.

Then we get by \eqref{eq3:maindelta} and by Theorem \ref{main} (2) below, that for $r$ sufficiently large the following holds\\

$r^{2-d}\cdot\bigl|\volu_{\Z} H_{r,M} - \volu H_{r,M}\bigr|\; \ll_d\; (b-a)\,\bar{q}^{d}q^{-1}r^{-\alpha K^{-1}}$ \beqa\phantom{.\hskip3mm} &+&  (b-a)\,\bar{q}^{d+1}q^{-1}\bigl(\abs{M}+2q^{-\frac{1}{2}}\frac{\abs{a}+\abs{b}}{r}\bigr)r^{-1}\; + \;\; T^{-\frac{1}{2}}\bar{q}^{\frac{d}{2}}q^{-1}\bigl(2 + q^{\frac{1}{2}}\bigr)^{d-2}\nn\\ &+& \; c(Q,M)\cdot r^{\alpha}\Bigl(1 + q_0^{-\frac{d}{2}}r^{\frac{d}{2}} + \bar{q}^dr^{\frac{d}{2}}(1+\log r )+ r^{d-2}\cdot\rho (r,Q,T)\Bigr)r^{2-d}\nn \eeqa
Taking the infimum over all $\alpha \in (0,1)$ we obtain\\

$r^{2-d}\cdot\bigl|\volu_{\Z} H_{r,M} - \volu H_{r,M}\bigr| \quad \ll_d\quad
 (b-a)\,\bar{q}^{d}q^{-1}r^{-\frac{1}{K}} $ \beqa & + &  (b-a)\,\bar{q}^{d+1}q^{-1}\bigl(\abs{M}+2q^{-\frac{1}{2}}\frac{\abs{a}+\abs{b}}{r}\bigr)r^{-1}
 +  T^{-\frac{1}{2}}\bar{q}^{\frac{d}{2}}q^{-1}\bigl(2 + q^{\frac{1}{2}}\bigr)^{d-2}\nn\\\phantom{.\hskip8mm} & +& \; c(Q,M)\cdot \Bigl(1 + q_0^{-\frac{d}{2}}r^{\frac{d}{2}} + \bar{q}^dr^{\frac{d}{2}}(1+\log r )+ r^{d-2}\cdot\rho (r,Q,T)\Bigr)r^{2-d}\nn \eeqa

By taking the infimum over all $T\geq 1$ we get with \eqref{defrho}\\

$r^{2-d}\cdot\bigl|\volu_{\Z} H_{r,M} - \volu H_{r,M}\bigr|
\quad \ll_d\quad (b-a)\,\bar{q}^{d}q^{-1}r^{-\frac{1}{K}} $\beqa +\; 
(b-a)\,\bar{q}^{d+1}q^{-1}\bigl(\abs{M}+2q^{-\frac{1}{2}}\frac{\abs{a}+\abs{b}}{r}\bigr)r^{-1} \, +
\, c(Q,M)\cdot\rho(r,\Q),\;\label{2.3main} \eeqa which proves 
Theorem \ref{maindelta} (2) for an appropriate choice of $r_0$.\bwend


The key tool in the previous proofs is the following
\begin{theorem}{\label{main}}
Let $\Q$ denote a non-degenerate $d$-dimensional quadratic form of
block-type, $d\geq 5$. Then for all $M\in\Rd$ there exist constants $c(\Q,M),
r_0 > 0$, such
that for any $r \geq r_0$ and any $T\geq 1$ \begin{itemize}
\item[(1)]
$ \Bigl|\volz{1}{\pm\eps} - \volr{1}{\pm\eps}\bigr|$\\ \phantom{.\hspace{20mm}.} $\ll_d\;
c(\Q, M)\cdot\eps^{-K}\Bigl( 1 + 2q_0^{-\frac{d}{2}}r^{\frac{d}{2}}+\bar{q}^{d+1}(\log q + 1)r^{d-2}\Bigl).$ 
\item[(2)] $ \Bigl|\volz{T^{-1/2}}{\pm\eps} - \volr{T^{-1/2}}{\pm\eps}\bigr|\;\ll_d \;  c(Q,M)\cdot\eps^{-K}$\\ \phantom{.\hspace{25mm}.} $\times \Bigl(1 + q_0^{-\frac{d}{2}}r^{\frac{d}{2}} + \bar{q}^dr^{\frac{d}{2}}(1+\log r )+ r^{d-2}\cdot\rho (r,Q,T)\Bigr)$,\\[2mm]
where $\rho (r,Q,T)$ is defined as in \eqref{defrhoT}.
\end{itemize}
\end{theorem}

\bw We want to estimate the difference between these two
approximations by integrals of theta functions. By
\eqref{g-integral}, \eqref{volzdef} and \eqref{volrdef} we have \vskip2mm

$ \Bigl|\volz{w}{\pm\eps} - \volr{w}{\pm\eps}\Bigr|\; =$\\
$$\Bigl|  \sum_{x\in \Z^d} \exp\Bigl[-\frac{2}{r^2}\,Q_+[x]\Bigr]\,\frac{1}{2\pi i} \int\limits^{\beta+ i\, \infty}_{\beta-i\, \infty}\exp\bigl[-z\cdot\Q[x-M]\bigr]\cdot h_{a,b,w}(z)\frac{dz}{z}\, \chi_{\pm\eps}\left(\frac{x}{r}\right)$$ $$ \hskip2mm -
 \int\limits_{\R^d} \exp\Bigl[-\frac{2}{r^2}\,Q_+[x]\Bigr]\,\frac{1}{2\pi i} \int\limits^{\beta+ i\, \infty}_{\beta-i\, \infty}\!\!\exp\bigl[-z\cdot\Q[x-M]\bigr]\cdot h_{a,b,w}(z)\frac{dz}{z}\, \chi_{\pm\eps}\left(\frac{x}{r}\right)dx\Bigr|.
$$
Choosing $\beta  = r^{-2}$, decomposing $Q[x-M] = Q[x]+Q[M]-2\langle x,QM\rangle$ (Recall, that $Q$ is self-adjoint.) and using  Fubini's theorem, we get\vskip2mm

$ \Bigl|\volz{w}{\pm\eps} - \volr{w}{\pm\eps}\Bigr|\;\; = \Bigl|\int\limits^{r^{-2}+ i\, \infty}_{r^{-2}-i\, \infty} \exp\bigl[-zQ[M]\bigr]h_{a,b,w}(z) $ $$\times \;\Bigl\{ \sum\limits_{x\in \Z^d} \exp\Bigl[-\frac{2}{r^2}\,Q_+[x] - zQ[x]+ i\langle x, 2t \Ima (z)QM\rangle\Bigr]\overline{\chi}_{\pm\eps}\left(\!\ffrac{x}{r}\!\right) \phantom{.....}$$
\beqa \hskip14mm
 - \int_{\Rd}\exp\Bigl[-\frac{2}{r^2}\,Q_+[x] - zQ[x]+ i\langle x, 2\Ima (z)QM\rangle\Bigr]\overline{\chi}_{\pm\eps}\left(\frac{x}{r}\right)dx\Bigr\}\frac{dz}{z}\Bigr|,\nn
\eeqa
where $\overline{\chi}_{\pm\eps}$ is defined as in \eqref{defbarchi}.\\ Since $ \overline{\chi}_{\pm\eps} (x)= \frac{1}{(2\pi)^d}\int_{\R^d}\widehat{\overline{\chi}}_{\pm\eps}(v)\exp[-i\langle x,v\rangle]dv$ holds by the Fourier inversion theorem, we obtain \\

$ \Bigl|\volz{w}{\pm\eps} - \volr{w}{\pm\eps}\Bigr| \;\; = $\\[4mm]
$ \Bigl|\int\limits^{r^{-2}+ i\, \infty}_{r^{-2}-i\, \infty}\exp\bigl[-zQ[M]\bigr]h_{a,b,w}(z) \frac{1}{(2\pi)^d}\;\;\int\limits_{\R^d} \widehat{\overline{\chi}}_{\pm\eps}(v)$
\beqa
&\times&\Bigl\{ \sum_{x\in \Z^d} \exp\Bigl[-\frac{2}{r^2}\,Q_+[x] - zQ[x]+ i\langle x, 2\Ima (z)\, QM-\frac{v}{r}\rangle\Bigr]\nonumber \\
&& \hskip3mm - \int_{\R^d}\exp\Bigl[-\frac{2}{r^2}\,Q_+[x] - zQ[x]+i\langle x, 2\Ima (z) QM-\frac{v}{r}\rangle\Bigr]dx\Bigr\}dv\;\frac{dz}{z}\Bigr|.\nn\\ &&\label{voluestimate1}
\eeqa

For $v\in \C^d$ we introduce the following theta sum and theta integral
\beqa
\theta_v(z) &\defi& \exp\bigl[-zQ[M]\bigr]\sum_{x\in\Z^d} \exp\left[- \Theta_{Q,M,r,v}(z,x)\right],\label{defthetasum} \\[2mm]
\theta_{0,v}(z) &\defi&\exp\bigl[-zQ[M]\bigr]\int_{R^d} \exp\left[- \Theta_{Q,M,r,v}(z,x)\right]dx\label{defthetaint}
\eeqa

where $\Theta_{Q,M,r,v}(z,x)\defi \frac{2}{r^2}Q_+[x] - z\cdot Q[x] - i\cdot\langle x,\frac{v}{r}- 2\Ima (z)QM\rangle$.\\
Then we can rewrite \eqref{voluestimate1} as follows\\[2mm]
$ \Bigl|\volz{w}{\pm\eps} - \volr{w}{\pm\eps}\bigr|$
$$= \Bigl|\int\limits^{r^{-2}+ i\, \infty}_{r^{-2}-i\, \infty}h_{a,b,w}(z) \frac{1}{(2\pi)^d}\int_{\R^d} \widehat{\overline{\chi}}_{\pm\eps}(v)\cdot \bigl\{\theta_v(z) - \theta_{0,v}(z) \bigr\}dv\;\frac{dz}{z}\Bigr|.
$$
Consider the segments $J_0\defi [r^{-2} - i\cdot\imbound ; r^{-2}+i\cdot\imbound] $ and $J_1 \defi \bigl(r^{-2} +i\,\R\bigr)\setminus J_0$. Then we may split\\[2mm]
$\Bigl|\volz{w}{\pm\eps} - \volr{w}{\pm\eps}\bigr|$ \beqa &\ll_d&\Bigl|\;\int\limits_{J_0}h_{a,b,w}(z) \frac{1}{(2\pi)^d}\int_{\R^d} \widehat{\overline{\chi}}_{\pm\eps}(v)\cdot \bigl\{\theta_v(z) - \theta_{0,v}(z) \bigr\}dv\;\frac{dz}{z}\nonumber\\ && - \int\limits_{J_1}h_{a,b,w}(z) \frac{1}{(2\pi)^d}\int_{\R^d} \widehat{\overline{\chi}}_{\pm\eps}(v)\cdot \theta_{0,v}(z) dv\;\frac{dz}{z}\nonumber \\ &&+\int\limits_{J_1}h_{a,b,w}(z) \frac{1}{(2\pi)^d}\int_{\R^d} \widehat{\overline{\chi}}_{\pm\eps}(v)\cdot \theta_v(z)dv\;\frac{dz}{z}\;\Bigr|\nn \\&=& \Bigl| I_0 - I_1 +I_2\Bigr|, \,\, \text{ say.}\label{eq:splitting}
\eeqa
Before estimating these integrals we derive a bound for $h_{a,b,w}(r^{-2} +it) , t \in \R$. Using
\beqa
\label{h-estimate1}
 \Bigl|\ffrac{\exp\{w\4 (r^{-2}+ i\4 t)\}-1}{w}\Bigr|
\leq \min \left \{ e\4\big|r^{-2} + i\4 t\big|, \ffrac{e+1}{w} \right \},
\eeqa
for \,$r^2 \geq \max (w, b) > 0$, \,$r \geq 1$, \,we obtain
\beqa \left|\frac{ h_{a,b,w}(r^{-2} +it)}{r^{-2}+it}\right| &\ll &\left( e^{br^{-2} }+e^{ar^{-2} } \right)\cdot \frac{1}{w\abs{r^{-2}+it}^2}\;\ll \;\;  \frac{1}{w\abs{r^{-2}+it}^2} ,\nn\\ &&\label{h-estimate2}\eeqa as well as \beqa
\left|\frac{h_{a,b,w}(r^{-2} +it)}{r^{-2}+it}\right| &\ll & \left( e^{br^{-2} }+e^{ar^{-2} }\right)\cdot\abs{r^{-2}+it}^{-1}\;\;\ll \;\; \abs{r^{-2}+it}^{-1}.\nn\\ &&\label{h-estimate3}
\eeqa
\vskip3mm
\noindent {\em Estimation of $ I_0$:} Inequality \eqref{h-estimate3} and Lemma \ref{thetaestimate1} for $t \in J_0$ yields
\beqa
\Theta_t &\defi& \left| \bigl(r^{-2}+i\,t\bigr)^{-1}h_{a,b,w}(r^{-2}+i\,t) \right|\nn \\ && \quad \quad\quad\times\quad\left|\int_{\R^d} \widehat{\overline{\chi}}_{\pm\eps}(v)\cdot \bigl\{\theta_v(r^{-2}+i\,t) - \theta_{0,v}(r^{-2}+i\,t) \bigr\}dv\right|\nn\\
&\ll_d& q_0^{-\frac{d}{2}}\abs{r^{-2}+it}^{-\frac{d+2}{2}}\exp\left[-c(Q)\cdot \Rea \bigl((r^{-2}+it)^{-1}\bigr)\right]\cdot\int_{\R^d} \left|\widehat{\overline{\chi}}_{\pm\eps}(v)\right|dv \nn\\ &&\hskip5mm + 2\cdot |r^{-2}+it|^{-1}\int_{\R^d}\left|\widehat{\overline{\chi}}_{\pm\eps}(v)\right| I_{(r,\infty)}(\abs{v})dv,\nn
\eeqa
 where $c(Q)$ is chosen according to Lemma \ref{thetaestimate1}.
Writing \,$\lvert \4r^{-2} + i\4t \4\rvert
= r^{-2}(1 +r^4\4t^2)^{1/2}$ and $ \Rea \bigl((r^{-2} +it)^{-1}\bigr) = \frac{r^2}{1+r^4t^2}$, \,we may
 introduce the variable $s = (1+ r^4\4t^2)^{-1}$ and the function
 \,$h(s)\defi s^{(d+2)/4}\4 \exp\{-c(Q)\4 s \4 r^2\}$. \,The maximal value of $h$ on $[\40, \infty)$ is attained at
\,${s_0= \ffrac{d+2}{4 \4 c(Q) \4 r^2}}$ and it is bounded by
\,$(c(Q)\4 r^2)^{-(d+2)/4}$ \, up to a constant depending on $d$ only.\\[2mm] Using the properties of $\chi_{\pm\eps}$ (see p. \pageref{taperingfunction}) and the fact that $\abs{v}\geq r$ implies $\abs{v}_{\infty}\geq d^{-1/2}r$ we now obtain
\beqa
\sup\limits_{t \in J_0} \Theta_t \!\!&\ll_d &\!\!q_0^{-\frac{d}{2}} \4 r^{d+2}\4
\sup\limits_{s\ge 0} \,h(s)\int_{\R^d} \left|\widehat{\overline{\chi}}_{\pm\eps}(v)\right|dv  + 2r^2\!\int_{\R^d}\left|\widehat{\overline{\chi}}_{\pm\eps}(v) \right|I_{(r,\infty)}(\abs{v})dv\nn\\
\!\! & \ll_d &\!\! q_0^{-\frac{d}{2}}\4 r^{d+2}\4 (c(Q) \4r^2 )^{-\frac{d+2}{4}}\cdot \!\int_{\R^d} \left|\widehat{\overline{\chi}}_{\pm\eps}(v)\right|dv \nn\\ &&\hskip40mm + \quad 2r^2\int_{\R^d}\left|\widehat{\overline{\chi}}_{\pm\eps}(v)\right| I_{(d^{-1/2}r,\infty)}(\abs{v}_{\infty})dv\nonumber\\
\!\!&\ll_d& \!\! q_0^{-d/2}\4 r^{d+2}\4 (c(Q) \4r^2 )^{- (d+2)/ 4}\cdot c_1(Q,M) \!\cdot\!\eps^{\epspot} + c_1(Q,M)\!\cdot\!\eps^{\epspot}\cdot r.\nn
\eeqa
Integrating this bound over $J_0$, we get for an appropriately chosen constant $c_2(Q,M)>0$
\beqa
|I_0| \leq \int\limits_{-\imbound}^{\imbound}\Theta_tdt \ll_d c_2(Q,M)\cdot\eps^{-K} q_0^{-\frac{d}{2}}r^{\frac{d}{2}} + c_1(Q,M)\cdot{\eps}^{\epspot}. \label{eq:I0estimate}
\eeqa

\noindent{\em Estimation of $ I_1$:}
Using Lemma \ref{thetasums2}, \eqref{thetaintegral} and \eqref{thetaestimate1-eq2}, we have
\beqa
\Bigl| \theta_{0,v}(z)\Bigr| \ll_d q_0^{-\frac{d}{2}}\abs{z}^{-\frac{d}{2}}.
\eeqa
Therefore, we get by the properties of $\chi_{\pm\eps}$ (see p. \pageref{taperingfunction}) and \eqref{h-estimate3} for $r^2\geq \max\{w,b,1\}$
\beqa
|I_1|&\ll_d& q_0^{-\frac{d}{2}}c_1(Q,M)\cdot \eps^{-K}\!\int_{J_1}\Bigl|\bigl(r^{-2}+i\,t\bigr)^{-\bigl( 1+\frac{d}{2}\bigr)}\Bigr| \, dt\nn\\ &\ll_d & q_0^{-\frac{d}{2}}c_1(Q,M)\cdot \eps^{-K}\int_{\imbound}^{\infty}t^{-\bigl(1+\frac{d}{2}\bigr)}dt\ll_d q_0^{-\frac{d}{2}}c_1(Q,M)\cdot \eps^{-K}r^{\frac{d}{2}},\nn\\[-3mm] &&\label{eq:I1estimate}
\eeqa
using the symmetry in t around 0.\\

\noindent {\em Estimation of $ I_2$:}
The estimate $ \bigl|h_{a,b,w}(r^{-2} +i\,t)\bigr|\ll_d \min\bigl\{1,\bigl(|r^{-2}+i\,t|w\bigr)^{-1}\bigl\}$
given by \eqref{h-estimate2} and \eqref{h-estimate3} implies
\beqa
|I_2| \!\!&\ll_d&\!\!\!\int_{\R^d}\!\int_{|t|>\imbound}\! \bigl| \theta_v\big(\frac{1}{r^2}+i\4t\big) \bigr|
\min \Bigl\{\,1, \ffrac{1}{w\,|\4 r^{-2}+i\4t\4|}\Bigr\}
\4 \ffrac{dt}{\abs{\4r^{-2}+i\4t\4}}  \left|\widehat{\overline{\chi}}_{\pm\eps}(v)\right|dv\nonumber
\\[2mm] \!\!&\ll_d&\!\!\!
\int_{\R^d}\!\int_{|u|>\imboundtwo} \bigl| \theta_v\bigl(r^{-2}+i
\4\pi\4 \frac{u}{2}\bigr)\bigr|\, g(u)\,du  \left|\widehat{\overline{\chi}}_{\pm\eps}(v)\right|dv, \label{I2-estimate1}
\eeqa
where
\beqa\label{g_min-def}
g(u) =\min \bigl\{1, (w\4|u|)^{-1}\bigr\}\, |u|^{-1}.
\eeqa
\vskip3mm
Using Lemma \ref{thetaestimate3} and the properties of $\chi_{\pm\eps}$ (see p. \pageref{taperingfunction}), we have
\beqa
|I_2|&\ll_d& {\bar{q}}^{d}\4 r^{d/2}\4 \int_{\R^d}\int_{|u|>\imboundtwo}(M_{1,t} \cdots M_{d,t})^{-1/2}\cdot g(u)\,du  \left|\widehat{\overline{\chi}}_{\pm\eps}(v)\right|dv\nn\\
&\ll_d& {\bar{q}}^{d}\4 r^{d/2}\cdot c_1(Q,M)\cdot\eps^{-K} \int_{|u|>\imboundtwo}(M_{1,t} \cdots M_{d,t})^{-1/2}\cdot g(u)\,du,\nn\\[-3mm] && \label{I2-estimate2}
\eeqa
where $M_{j,t}$ denote Minkowski's successive minima for the norm on  ~$\mathbb{R}^{2d}$  related to~$\Q$, defined
by \eqref{eq:norm} and \eqref{eq:forms} and $c_1(Q,M)>0$ is  a constant chosen according to  \eqref{defbarchi}.
Denote
\beqa
G(\kappa,\xi)\defi \int^{\xi}_{\kappa} g(t)\,dt, \qquad\text{for}\; 0<\kappa<\xi\le\infty.
\eeqa
For  $\kappa \ge \xi>0$ we define $G(\kappa,\xi) =0$.
Note that
\beqa\label{eq:G}
G(\kappa,\xi) = \left\{ \begin{array}{l@{\;  \;}l}
\;\; \log (\xi/\kappa), & \text{for } \,\, \kappa \leq \xi \leq w^{-1}, \\[0.2cm]
- \log (w \4\kappa) + 1-(w\4\xi)^{-1}, & \text{for } \,\, \kappa\le w^{-1}\le \xi,\\[0.2cm]
(w\4\kappa)^{-1}-(w\4\xi)^{-1}, & \text{for } \,\,w^{-1}\le  \kappa\le \xi.
\end{array} \right.
\eeqa
The equality \eqref{eq:G} and the definition of the function
\,$G$ \,imply the bound
\beqa\label{eq:GGG}
G(\kappa, \xi)\le \min\bigl\{\bigl|\log (w \4\kappa)\bigr| + 1,\
\bigl|\log (\xi/\kappa)\bigr|,\ (w\4\kappa)^{-1}\bigr\}
\quad\text{for}\quad \kappa, \xi >0.\nn\\[-1mm]
\eeqa
Writing  $$ {M}(t) = M_{1,t} \cdots M_{d,t},$$
the upper bound for $|I_2|$ in \eqref{I2-estimate2}
in terms of Minkowski's successive minima now
yields
\beqa
\label{eq:2.22}
| I_2 |&\ll_d &{\bar{q}}^d\4r^{d/2}\cdot c_1(Q,M)\cdot \eps^{-K} \int_{|t|>\imboundtwo}
\ffrac{g(t)}{{M}(t)^{1/2}}\4 dt\nn\\
&=&2\4 {\bar{q}}^d\4r^{d/2}\cdot c_1(Q,M)\cdot\eps^{-K} I_3,\quad
\eeqa
where
\beqa
\label{eq:2.22b}
 I_3  =\int_{\imboundtwo}^\infty
\ffrac{g(t)}{{M}(t)^{1/2}}\4 dt.
\eeqa
The last equality in \eqref{eq:2.22} follows from the fact that
the functions \,$g(\cdot )$ \,and  \,${M}(\cdot)$ \,are even
(see \eqref{eq:symm}).\\[2mm]
After this preparations, we may now complete the proof of Theorem \ref{main}:\\

\noindent \textbf{Proof of  Theorem \ref{main} (1)}.\quad\\
Let 
\beqa
\gamma({\kappa,\xi}) 
= r^d\4 \inf_{\kappa \le t\le \xi} {M}(t), \qquad\text{for}\; \kappa , \xi \in\R.
\label{eq:gamM}
\eeqa

Applying Lemma \ref{final} for the interval with endpoints ~\,$\kappa= \imboundtwo$ and 
 ~$\xi=\infty $, \,we get
\beqa
I_3
&\ll_d& q_0^{-1} r^{d/2-2} \4 \int_{\gamma_0}^{D_0}v^{-1/2 + 1/d}
\bigl(q {v}^{1/d}\4 G\bigl(\kappa_0(v^{1/d}),\infty\big)
+ g\bigl(\kappa_0(v^{1/d})\bigr)\bigr) 
 \ffrac {dv} v \nn\\ && \hskip6cm + G(\imboundtwo ,\infty)\label{eq:Intg}
\eeqa
 \vskip-5mm with  
\beqa
\gamma_0=\gamma({\imboundtwo ,\,\infty}),\;
D_0=\max\bigl\{\left(\frac{r}{2d}\right)^d,\gamma_0\bigr\},
\; \kappa_0(v) = \max\bigl\{ \imboundtwo , \frac{1}{2\4q\4  v\4d^{1/2}}\bigr\}.\nn\\ \label{eq:gDa}
\eeqa
Note that \,$\gamma_0 \geq 1$ \,by
\eqref{eq:minima0}.
 In the sequel we choose  $w=1$. Using \eqref{g_min-def}, \eqref{eq:GGG}, 
\eqref{eq:2.22}, \eqref{eq:Intg}, \eqref{eq:gDa} and hence $g\bigl(\kappa_0(v^{1/d})\bigr) \ll_d q\, v^{1/d}$, we obtain for $d>4$ and $ r\geq \max\left\{\frac{2}{\pi\,q}; \frac{2}{\pi}\right\}:$
\beqa
|I_2 |&\ll_d& c_1(Q,M)\cdot\eps^{-K} {\bar{q}}^{d+1}\4r^{d/2}\4 r^{d/2-2} \4 
\int_1^{D_0}v^{-1/2 + 2/d}\4 \bigl(\log (q\4v^{1/d}) + 2\bigr) 
 \ffrac {dv} v \nn\\ && \phantom{.\hskip5cm} + c_1(Q,M)\cdot\eps^{-K}{\bar{q}}^d\4r^{d/2}\4(\log r+1)\nonumber \\[2mm]
& \ll_d & c_1(Q,M)\cdot\eps^{-K}{\bar{q}}^{d+1}\4(\log q+1)\, r^{d-2}.
\eeqa 
For  $r\geq r_0\defi \max\left\{\frac{2}{\pi\,q}; \frac{2}{\pi}, r_0(Q,M)\right\}$, where $r_0(Q,M)$ is a constant chosen as in Lemma \ref{thetaestimate1} and \ref{thetaestimate2}, this bound
for $I_2$ yields in view of \eqref{eq:splitting}, \eqref{eq:I0estimate} and \eqref{eq:I1estimate}, that
\beqa  \Bigl|\volz{1}{\pm\eps} - \volr{1}{\pm\eps}\bigr|\;\ll_d\; c_2(Q,M)\cdot\eps^{-K}q_0^{-\frac{d}{2}}r^{\frac{d}{2}}\hskip20mm \nn\\ +\; c_1(Q,M)\cdot\eps^{-K}\bigl( 1 + q_0^{-\frac{d}{2}}r^{\frac{d}{2}} +  {\bar{q}}^{d+1}\4(\log q+1)\, r^{d-2}\bigr),\nn
\eeqa
where the constants $c_1(Q,M)$ and $c_2(\Q,M)$ are chosen according to Lemma \ref{defbarchi} and \eqref{eq:I0estimate}. Setting $c(Q,M)\defi \max\{c_1(Q,M),c_2(Q,M)\}$,
this proves Theorem \ref{main}  (1).
\medskip

\noindent \textbf{Proof of Theorem \ref{main} (2)}.\\
In order to use nontrivial bounds for $\gamma(\kappa,\xi)$ in the irrational
case we introduce further auxiliary parameters $\eta, T$ 
such that $\imboundtwo  \le \eta\le T$ \,with \, $T\ge1$ \,which will be 
determined and optimized later. 
Thus we may
split the integral $I_3$ in \eqref{eq:2.22b}
which bounds  $|I_2|$ in \eqref{eq:2.22} into the parts
\beqa
I_3 & = & \int_{\imboundtwo}^\eta \ffrac{g(t)}{{M}(t)^{1/2}} \, dt 
+\int_\eta^T \ffrac{g(t)}{{M}(t)^{1/2}} \, dt + \int_T^\infty  
\ffrac{g(t)}{{M}(t)^{1/2}} \, dt\nonumber \\[2mm] 
 & = & I_4 + I_5 + I_6, \qquad \text{say}.\label{eq:I456}
\eeqa

We define similarly to \eqref{eq:gDa}
\beqa
\gamma_1=\gamma({\imboundtwo ,\eta}),\;
\gamma_2=\gamma({\eta,T}),\;
\gamma_3=\gamma({T,\infty}),\;\label{eq:gammas}
\eeqa
\beqa
 D_j=\max\bigl\{(2d)^{-d}\4 r^d,\,\gamma_j\bigr\},
\qquad j=1,2,3,
\label{eq:MJ}
\eeqa
$$
\kappa_1(v)=\max\bigl\{ \imboundtwo , f(v)\bigr\},\;
\kappa_2(v)=\max\bigl\{ \eta, f(v)\bigr\},\;
\kappa_3(v)=\max\bigl\{ T, f(v)\bigr\},\nn\\[-1mm]
$$\vskip-7mm
\beqa
\label{eq:av}
\eeqa
where
\,$f(v)= (2\4q\4  v\4d^{1/2})^{-1}$, \,$v>0$.
By \eqref{eq:minima0} we have again
\beqa
\label{eq:gammabound}
\gamma_j\ge1,\qquad j=1,2,3.
\eeqa
Using \eqref{g_min-def} and  \eqref{eq:av}, we see that
\beqa
\label{eq:gbound}
g(\kappa_j(v))\le 2\4q\4  v\4d^{1/2} ,\qquad j=1,2,3.
\eeqa

First, we apply Lemma \ref{final} as above to the interval 
with endpoints $\kappa=\imboundtwo $ and $\xi=\eta$. Corollary 
\ref{multineq} implies that, if $\eta \geq \ffrac{\pi}{2\,q_0\,q\,r}$
the quantity \,$\gamma_1$ \,(defined by \eqref{eq:gamM} and
\eqref{eq:gammas}) satisfies
\beqa\label{eq:gb}
\gamma_1 \ge \delta\defi (d\4q\4\eta)^{-d},
\eeqa
since \,$d\ge 5$ and $\inf_{t\in [\frac{2}{\pi\, r},\eta]}\left\{\frac{q_0|t|r}{2};\frac{1}{q|t|r}\right\} = \ffrac{1}{q\,\eta\,r} $, whenever $\eta \geq \ffrac{\pi}{2\,q_0\,q\,r}$.

Lemma \ref{final}
 yields in view of \eqref{eq:G}, \eqref{eq:GGG}, \eqref{eq:gbound}
and \eqref{eq:gb}  the estimate
\beqa
I_4 &\ll_d&  q_0^{-1}r^{d/2-2} \4 \int_{\gamma_1}^{D_1}v^{-1/2 + 1/d}
\bigl( {v}^{1/d}\4q\, G\bigl(\kappa_1(v^{1/d}),\eta\bigr)
+ g\bigl(\kappa_1(v^{1/d})\bigr)\bigr) 
 \ffrac {dv} v \nn\\ && \hskip2cm+G(\imboundtwo ,\eta) \nonumber\\
& \ll_d& q_0^{-1}\4 q\4r^{d/2-2} \4 \int_{\delta}^{D_1}v^{-1/2 + 2/d}\4 
\bigl(\bigl|\log (q\4v^{1/d}\4\eta)\bigr| + 1\bigr)\4 
\ffrac {dv} v + G(\imboundtwo ,\eta)\nonumber \\ 
\label{eq:I3}
& \ll_d& q_0^{-1}q^{d/2-1}\4r^{d/2-2}\4\eta^{d/2-2}+G(\imboundtwo ,\eta),
\eeqa
 provided that $d>4$, using the change of variables 
$ v =\delta \4 u$ in the last inequality.

In order to estimate $I_5$ we choose $\kappa=\eta$, and $\xi=T$. By Lemma \ref{final} we obtain as above
\beqa
I_5 &\ll_d&  q_0^{-1}r^{d/2-2} \4 \int_{\gamma_2}^{D_2}v^{-1/2 + 1/d}
\bigl( {v}^{1/d}\4q\, G\bigl(\kappa_2(v^{1/d}), T)+ g\bigl(\kappa_2(v^{1/d})\bigr)\bigr) 
 \ffrac {dv} v  \nn\\ && \hskip2cm +\; G(\eta,T) \nonumber\\
& \ll_d& q_0^{-1}q\4r^{d/2-2}\int_{\gamma_2}^{D_2}v^{-1/2 + 2/d}\4\bigl(
\bigl| \log (q\4v^{1/d}/w)\bigr| + 1\bigr)\ffrac {dv} v 
+G(\eta,T)\nonumber \\
\label{eq:I4}
& \ll_d&  q_0^{-1}q\,r^{d/2-2}\4 {\gamma_2}^{-1/2 +2/d}\4 
\bigl(\bigl|\log(q\4\gamma_2)\bigr| +\bigl|\log w\bigr| +1\bigr)+ G(\eta,T).
\eeqa

Finally for the term ~$I_6$ choose $\kappa = T$ 
and $\xi = \infty$~ and use \eqref{eq:gammabound}
for $j=3$. 
Recall that we choose $T\ge 1$. Thus, similarly 
as above, using Lemma \ref{final} and the fact, that
$G(\kappa_3(v^{1/d}), \infty)\le G(T, \infty)\le T^{-1}\4w^{-1}$ 
and $ g\big(\kappa_3(v^{1/d})\big)\le T^{-2}\4 w^{-1}$, \,we obtain
(see  \eqref{g_min-def}, \eqref{eq:GGG} and \eqref{eq:av})
\beqa
I_6 &\ll_d&  q_0^{-1} r^{d/2-2} \4 \int_{1}^{D_3}v^{-1/2 + 1/d}
\bigl( {v}^{1/d}\4q\, G\bigl(\kappa_3(v^{1/d}), \infty)
+ g\bigl(\kappa_3(v^{1/d})\bigr)\bigr) 
 \ffrac {dv} v \nn\\ && \hskip2cm +G(T,\infty) \nonumber\\
\label{eq:I5}
&\ll_d&  q_0^{-1}q\4 r^{d/2-2}\4T^{-1}\4w^{-1} + G(T,\infty).
\eeqa

Collecting \eqref{eq:I3}--\eqref{eq:I5},
 we get by combining the terms $G(\kappa,\xi )$ and using \eqref{eq:I456}
and the estimates \eqref{eq:gammabound}
\beqa
I_3 &\ll_d& q_0^{-1}r^{\frac{d}{2}-2}\4 \Bigl\{
 q^{\frac{d}{2}-1}\4 \eta^{\frac{d}{2}-2}+ q\4 \gamma_2^{-\frac{1}{2}+\frac{2}{d}}\4
\bigl(\log(q\4\gamma_2) +|\log w\4| +1\bigr)+ 
\frac{q}{Tw} \Bigr\}\nonumber\\ &&\hskip1cm + \; G(\imboundtwo , \infty).\label{eq:I2}
\eeqa

In view of \eqref{eq:2.22} this bound for $I_3$ yields 
\beqa 
|I_2|&\ll_d& c_1(Q,M)\cdot\eps^{-K}\cdot q^{d}\4r^{\frac{d}{2}}\bigl(1+\log r \bigr) \; + \; c_1(Q,M)\cdot\eps^{-K}\cdot q_0^{-1}\bar{q}^d\cdot r^{d-2}\nn\\&& \quad \times \;\Bigl\{(T\4 w)^{-1} + q^{\frac{d}{2}-1}\4 \eta^{\frac{d}{2}-2}
+ q\4{\gamma_2}^{-\frac{1}{2}+ \frac{2}{d}}\4
\bigl(\log(q\4\gamma_2) +|\log w\4| +1\bigr) \Bigr\} \nonumber\\
&\ll_d& c(Q,M)\cdot \eps^{-K}\cdot \bar{q}^{d}\4r^{\frac{d}{2}}\bigl(1+\log r \bigr)\; + \; c(Q,M)\cdot \eps^{-K}\cdot r^{d-2}\nn\\ && \quad\quad \times \Bigl\{\frac{\bar{q}^{d+1}}{T\4 w} + \bar{q}^{\frac{3d}{2}}\4 \eta^{\frac{d}{2}-2}
+ \bar{q}^{d+2}\4{\gamma_2}^{-\frac{1}{2}+ \frac{2}{d}}\4
\bigl(\log(\bar{q}\4\gamma_2) +|\log w\4| +1\bigr) \Bigr\}
,\nn\\ &&
\label{eq:finalraw}
\eeqa
where $c(Q,M) \defi \max\{c_1(Q,M),c_2(Q,M)\}$.
By Lemma \ref{irrational} for $\eta, T$ fixed, 
we have $\gamma_2 \rightarrow \infty$ for 
 $ r\to \infty$ and we may now choose the auxiliary
 parameters $ \eta, w$ and $T$ to minimize the right hand side of
\eqref{eq:finalraw} as follows. 
Let 
\beqa
T \geq 1,\; w = T^{-1/2}, \; 
\eta = \max\bigl\{\imboundtwo ,\,\frac{\pi}{2q_0\,q\,r},\,T^{-\frac{1}{d-4}}\bigr\},\; 
\label{eq:finall}
\eeqa 
provided that $d \ge5$.\\
For $r\geq r_0\defi \max\Bigl\{ \frac{2}{\pi} ,\,\frac{\pi}{2q_0\,q}, r_0(Q,M)\Bigr\}$, where $r_0(Q,M)$ is a constant chosen as in Lemma \ref{thetaestimate1} and \ref{thetaestimate2}, we obtain in view of \eqref{eq:splitting}, \eqref{eq:I0estimate}, \eqref{eq:I1estimate}, \eqref{eq:gammas}, \eqref{eq:gammabound},
 \eqref{eq:finalraw} and \eqref{eq:finall}  the following bound:\\[2mm]
$  \Bigl|\volz{T^{-1/2}}{\eps} - \volr{T^{-1/2}}{\eps}\Bigr|$\\
\phantom{.\hskip12mm .}$\ll_d \;\; c(Q,M)\cdot\eps^{-K}\Bigl(1 + q_0^{-\frac{d}{2}}r^{\frac{d}{2}} + \bar{q}^dr^{\frac{d}{2}}(1+\log r )+ r^{d-2}\cdot\rho (r,Q,T)\Bigr)$,\\[2mm]
where $\rho (r,Q,T)$ is defined as in \eqref{defrhoT}.
This completes the proof of Theorem \ref{main}  (2).\bwend

\bwthm{Theorem \ref{diophant}}\\ The estimate \eqref{eq:metric1}
immediately follows from Corollary \ref{lemmadiophantine1}. This
inequality ensures that there exists a $t\in [\kappa,\xi]$ such that
~$M_{1,t} > \tau$ whenever $c(d)\bigl(q\4
\tau^2\4(\xi-\kappa)+\tau\4 r^{-1}\bigr) < \xi -\kappa$. This
condition is equivalent to \beqa\tau < \left(\frac{1}{c(d)} -
q\4\tau^2\right)\4(\xi-\kappa)\4r.\label{eq:diophant1}\eeqa Due to
the fact, that $\tau \leq \tau_Q$, where $ \tau_Q\defi
\left(\frac{c(d)+2}{2c(d)q}\right)^{\frac{1}{2}}$, implies
$\frac{1}{c(d)} - q\4\tau^2 \ge \frac{1}{2}$, we may conclude, that the
condition \eqref{eq:diophant1} (and hence $M_{1,t} > \tau$) follows
from the inequality
 $\tau\leq \min\left\{\tau_Q , r(\xi-\kappa)/2\right\},$ which proves
\eqref{eq:metric2}.\\
By definition of $M_{1,t}$ the inequality
$M_{1,t}>\overline{\tau}\defi\min\{\tau_Q, r\4(\xi-\kappa)/2\}$
implies that if
 $0 <|n|_{\infty} < \overline{\tau}\4 r$ then $\overline{\tau}r\4 \| t\4 \Q\4 n \| >\overline{\tau}^2$.
For $\nu > \tau_Q$ exists a $ r\geq 1$ such that $\nu =
\overline{\tau}\4 r$. Therefore, we get by \eqref{defD(t,v)} that $
D(t,v) \geq \overline{\tau}^2$. Furthermore, we have
$\overline{\tau}^2= \min\{ \tau_Q^2, \nu \4 (\xi - \kappa)/2\},$
since either $r\4(\xi-\kappa)/2 > \tau_Q$ and $\overline{\tau} =
\tau_Q$ or $\overline{\tau} = r\4(\xi-\kappa)/2$ otherwise. This
proves \eqref{eq:metric3}. \phantom{ffgfgfggfg}\bwend

\bwthm{Theorem \ref{thdistr}}\\
Since the cube $C_r$ is compact the quantity \beqa a_r \defi \min\big\{\Q[x-M]: x \in C_r \big\}\label{defar}\eeqa
is a well-defined real number and we obviously get
\beqa F_{r,M}(b) = H_{r,M}^{a_r,b},\label{defFrMb}
\eeqa
where $H_{r,M}^{a_r,b}$ is defined as in \eqref{defHrab}. 

A careful analysis of the proof shows, that Theorem \ref{main} also holds for $a = a_r, r\geq r_0$. This, together with Lemma \ref{estmaindistr} yields that for $K=K(d)$ chosen according to \eqref{taperingfunction} there exist
 constants \,$c_j>0$, \,$j=1,...,5$, \,depending on $Q$ and $d$ only
and a constant $r_0 = r(\Q, M,b)>0$ such that, for any $r \ge r_0,$ it holds (cf. proof of Theorem \ref{maindelta}):
\begin{itemize}
\item[(1)]$ \bigl|\volu_{\mathbb{Z}}F_{r,M}(b)-\volu F_{r,M}(b)\bigr|\;\le \; r^{d-2}\cdot \bigl(c_1\cdot(b-a_r+1) + c_2\bigr).$\\[-1mm]

\item[(2)]$
\bigl|\volu_{\mathbb{Z}}F_{r,M}(b)-\volu F_{r,M}(b)\bigr| $ $$\;\;\le
\;\; r^{d-2} \cdot \bigl(
c_3\cdot (b-a_r)\,r^{-\frac{1}{K}} + 
c_4\cdot (b-a_r)\,r^{-1} \, +
\, c_5\cdot\rho(r,\Q )\bigr),$$ where $ \lim\limits_{r\to \infty}
\rho(r, \Q) =0$,\;\; provided that\;$  \Q $ is irrational.
\end{itemize}
Dividing these inequalites by the inequality in Lemma \ref{volestdistr} (2) for $\xi =1$ completes the proof of Theorem  \ref{thdistr}. \bwend

\section{Lemmas}
In the sequel, let $I = [a,b], a,b \in\R$ and $I_0$ denote finite
intervals. For $M \in \Rd$ we consider \beqa H(r)\defi H(r, I_0,I,M
)\defi \bigl\{x \in \Rd \, : \, r^{-1}\abs{x}_{\infty} \in I_0 ,\;
Q[x-M] \in I\bigr.\bigr\}.\label{defHr} \eeqa

The diagonal matrix $D(Q)$ is defined by
$$ \bigl( D(Q) \bigr)_{i,j}\defi \left\{\begin{array}{ll} \sqrt{\abs{q_i}},& \text{ if } j=i,\\0,& \text{ otherwise},\end{array}\right.\;\;\; 1\leq i,j\leq d.$$

\begin{lemma}\label{vol3}
Let $I_0 = [0,\xi ]$ and $\tau = \xi + \frac{\abs{D(Q)M}}{r}, \sigma = q_0^{d/2}\xi - \frac{\abs{D(Q)M}}{r}$. For the volume of $H(r)$ defined in \eqref{defHr} it holds
\beqa \volu \, H(r) \ll_d (b-a)q_0^{-d/2}q^{(d-2)/2}\tau^{d-2}r^{d-2}.\nn
\eeqa
 If $\sigma>0$ and $\abs{a}+\abs{b}\leq \sigma^2r^2/5$ then
 \beqa \volu\, H(r)\gg_d (b-a)q^{-d/2}\sigma^{d-2}r^{d-2}.\nn\eeqa
\end{lemma}
\bw \cite{bentkus-goetze:1999}, p. 1023, Lemma 8.2 or \cite{elsner:2006}, p. 24, Lemma 2.4.3 \bwend

\begin{lemma}\label{vol4}
Let $I_0  = [1-\delta, 1+\delta], \; 0\leq \delta\leq 1/4$. Assume that $r$ is large enough, that
\beqa \eps_1\defi r^{-1}|D(Q)M|\leq q_0^{1/2}/4\;\;\; \text{ and }\;\;\;\eps_2\defi  r^{-2}\left( |a|+|b|\right)\leq \frac{1}{8}q_0\label{epsbounds}
\eeqa
holds. Then for the volume of $H(r)$ defined in \eqref{defHr} it holds
\beqa \volu \, H(r) \ll_d (b-a)\left(\delta + q_0^{-1/2}\eps_1 + 2q_0^{-1/2}\eps_2\right)\, r^{d-2}q_0^{-d/2}q^{(d-2)/2}.\nn
\eeqa
\end{lemma}
\bw \cite{bentkus-goetze:1999}, p. 1025, Lemma 8.3 or \cite{elsner:2006}, p. 26, Lemma 2.4.4 \bwend

Due to the fact, that for $a_r$ defined as in \ref{defar} the inequality
\beqa \ffrac{|a_r|}{r^2} \leq q \eeqa
holds for $r$ large enough, we obtain in the case $a = a_r$  the following lemma by slightly modifying the proof of Lemma \ref{vol3} given in \cite{bentkus-goetze:1999} resp. \cite{elsner:2006}. Using these modifications we also get an analog result as in Lemma \ref{vol4}. 
\begin{lemma}\label{volestdistr} Let $I_r\defi [a_r,b]$. There exist constants $C_{Q,1}, C_{Q,2}\geq 1$ depending on $d$  and $Q$ only and a constant $r_0 = r_0(Q,M,b)\geq 1$ such that for $r\geq r_0$ the volume of  $F(r) \defi H(r,I_0,I,M)$  defined as in \eqref{defHr} can be estimated as follows:
\begin{itemize}
\item[(1)] For $I_0 = [0,\xi]$ it holds: $\volu \, F(r) \leq (b-a_r)\cdot C_{Q,1}\cdot \xi^{d-2}r^{d-2}$.\vskip2mm
\item[(2)] For $I_0 = [0,\xi]$ it holds: $ \volu\, F(r)\geq (b-a_r)C_{Q,1}\cdot\xi^{d-2}r^{d-2}$\vskip2mm
\item[(3)] For $I_0 = [(1-\delta), 1+\delta)], 0\leq \delta\leq 1/4$ it holds: 
$$\volu \, F(r) \leq (b-a_r)C_{Q,2}\cdot\delta \cdot r^{d-2}.$$
\end{itemize}
The constants  $C_{Q,1}, C_{Q,2}$ can be computed explicitly.
\end{lemma}
\vskip3mm


In the sequel we want to estimate the error terms caused by the approximations of the (lattice point) volumes of the hyperbolic shell $H_{r,M}$:\\

In the notation of \eqref{volzdef}-\eqref{volrdef}, considering for $\eps>0$
$$\psi_{r,\pm\eps}(x) =  \exp\bigl[-\frac{2}{r^2}Q_+[x]\bigr]\chi_{\pm\eps}\bigl(\frac{x}{r}\bigr)$$ and $$
\Delta_{\pm\eps} = \Bigl|\int_{\Rd} I_{H_{r,M}}(x)\psi_{r,\pm\eps}(x)dx - \sum_{x\in\Zd} I_{H_{r,M}}(x)\psi_{r,\pm\eps}(x)\Bigr|
,$$ defined as in \eqref{defpsi_r} and \eqref{defDeltaeps}, respectively, we define additionally \beqa
v_{\eps}&\defi& \volu \Bigl( H_{r,M} \cap \bigl\{x\in \Rd \bigl| r(1-\eps)\leq \abs{x}_{\infty}\leq r(1+\eps)\bigr.\bigr\}\Bigr)
\eeqa
and get the following estimate
\begin{lemma}\label{approx1} For $0<\eps \leq \frac{1}{4}$ there exists a constant $c = c(d) >0$ such that\\[3mm]
$\bigl|\volu_{\Z} H_{r,M} - \volu H_{r,M}\bigr| \;\;\leq \;\;\max \bigl\{ \Delta_{-\eps} \,; \,\Delta_{\eps}\bigr\} \; +\; c\cdot (b-a)q_0^{-d/2}q^{(d-2)/2}$ \beqa \hskip2cm\times \left(\eps + q_0^{-1/2}q^{1/2}|M|r^{-1} + 2q_0^{-1/2}(|a|+|b|)r^{-2}
\right)\, r^{d-2}.
\label{approx1claim}
\eeqa
\end{lemma}
\bw Obviously, we can estimate
\beqa
 \volu_{\Z} H_{r,M}\!\! &\leq& \!\!\sum_{x\in\Zd} I_{H_{r,M}}(x)\psi_{r,\eps}(x),\;\;\;\;\volu H_{r,M}\;\leq\;\int\limits_{\Rd} \!I_{H_{r,M}}(x)\psi_{r,-\eps}(x)dx  + v_{\eps},
\nn\\
\volu_{\Z} H_{r,M} \!\!&\geq &\!\!\sum_{x\in\Zd} I_{H_{r,M}}(x)\psi_{r,-\eps}(x),\;\; \;
\volu H_{r,M}\geq \int\limits_{\Rd} \!I_{H_{r,M}}(x)\psi_{r,\eps}(x)dx  - v_{\eps}.\nn\eeqa
If $\volu_{\Z} H_{r,M} - \volu H_{r,M}\geq 0$ these estimates imply
\beqa
\bigl|\volu_{\Z} H_{r,M} - \volu H_{r,M}\bigr|\leq \Delta_{+\eps} + v_{\eps},\nn
\eeqa
and otherwise we obtain
\beqa
\bigl|\volu_{\Z} H_{r,M} - \volu H_{r,M}\bigr| \leq \Delta_{-\eps} + v_{\eps}.\nn
\eeqa
Using Lemma \ref{vol4} for $I_0 = [1-\eps,1+\eps]$ we get since $|D(Q)M|\leq q^{1/2}|M|$ that
\beqa
v_{\eps} \ll_d (b-a)\left(\eps + q_0^{-1/2}q^{1/2}|M|r^{-1} + 2q_0^{-1/2}(|a|+|b|)r^{-2}
\right)\, r^{d-2}q_0^{-d/2}q^{(d-2)/2},\nn
\eeqa
which proves \eqref{approx1claim}.
\bwend
\begin{lemma}\label{approx2}
For fixed $a,b \in \R, w>0$ and the functions $g$ defined in
\eqref{defg} the following holds
\begin{itemize}
\item[(1)] There exist $ a' \in [a-w;a+w]$ and $b' \in [b-w;b+w]$ such that $$\sum\limits_{x \in \Zd}\bigl(I_{[a,b]} - g_{a',b', w}\bigr)\bigl(Q[x-M]\bigr)\psi_{r,\pm\eps}(x) = 0.$$
\item[(2)]$\hskip-4mm\sup\limits_{{a' \in [a-w;a+w]}\atop{b' \in [b-w;b+w]}}\!\Bigl|\int\limits_{\Rd}\bigl(I_{[a,b]} - g_{a',b', w}\bigr)\bigl(Q[x-M]\bigr)\psi_{r,\pm\eps}(x)dx\Bigr| $ \vskip-3mm$$\hskip3cm\ll_d  8w q_0^{-\frac{d}{2}}q^{\frac{d-2}{2}}\bigl(1+\eps+q^{\frac{1}{2}}\frac{\abs{M}}{r}\bigr)^{d-2}r^{d-2}\! .$$
\end{itemize}
\end{lemma}
\bw The sum in (1) is finite, since $\psi_{r,\pm\eps}$ has bounded support. Hence, the map
$$ G: (a',b') \longmapsto \sum\limits_{x \in \Zd}\bigl(I_{[a,b]} - g_{a',b', w}\bigr)\bigl(Q[x-M]\bigr)\psi_{r,\pm\eps}(x)
$$ is continuous and (1) follows by the intermediate value theorem.\\
For all $ a' \in [a-w;a+w]$ and all $b' \in [b-w;b+w]$ we can estimate
\beqa
\Bigr|(I_{[a,b]} - g_{a',b', w}\bigr)\bigl(Q[x-M]\bigr)\Bigr| \leq I_{\bigl([a-2w; a+2w]\;\cup \;[b-2w ;  b+2w]\bigr)}\bigl(Q[x-M]\bigr).
 \eeqa
This implies\\[2mm]
$\sup\limits_{a' , b'}\left|\int\limits_{\Rd}\bigl(I_{[a,b]} - g_{a',b', w}\bigr)\bigl(Q[x-M]\bigr)\psi_{r,\pm\eps}(x)dx\right| $
\beqa&\leq &\int  I_{\bigl([a-2w; a+2w]\;\cup \;[b-2w ;  b+2w]\bigr)}\bigl(Q[x-M]\bigr)\psi_{r,\pm\eps}(x)dx\nn\hskip15mm\\ &\leq& \int  \bigl(I_{[a-2w; a+2w]}+ I_{[b-2w ;  b+2w]}\bigr)\bigl(Q[x-M]\bigr)I_{\bigl[0 ; r(1+\eps)\bigr]}(\abs{x}_{\infty})dx,\hskip15mm
\label{supest}\eeqa
since $ \psi_{r,\pm\eps}(x)\leq I_{\bigl[0 ; r(1+\eps)\bigr]}(\abs{x}_{\infty})$. \\
Using Lemma \ref{vol3} with $I_0 = [0,1+\eps]$,  we get by \eqref{supest}\\

$ \sup\limits_{a', b'} \Bigl|\int\limits_{\Rd}\bigl(I_{[a,b]} - g_{a',b', w}\bigr)\bigl(Q[x-M]\bigr)\psi_{r,\pm\eps}(x)dx\Bigr|$
\vskip-4mm\beqa\hskip1cm&\ll_d& 8w q_0^{-d/2}q^{(d-2)/2}\bigl(1+\eps+r^{-1}\abs{D(Q)M}\bigr)^{d-2}r^{d-2}\nn\\[2mm]&\leq& 8w q_0^{-d/2}q^{(d-2)/2}\bigl(1+\eps+r^{-1}q^{1/2}\abs{M}\bigr)^{d-2}r^{d-2},
\eeqa
which proves (2).
\bwend\vskip-10mm \phantom{....}
\begin{lemma}\label{approx3}
Consider $\Delta_{\pm\eps} , \eps>0$  defined in \eqref{defDeltaeps}. Then:\\[2mm]
$
\Delta_{\pm\eps} \ll_d \sup\limits_{{a' \in [a-w;a+w]}\atop{b' \in [b-w;b+w]}}\!\Bigr|\vrapprox{w}{\pm\eps}{a'}{b'} - \vzapprox{w}{\pm\eps}{a'}{b'}\Bigr| $\\[-2mm] $\phantom{...}\hfill + \; 8w q_0^{-\frac{d}{2}}q^{\frac{d-2}{2}}\bigl(1+\eps+q^{\frac{1}{2}}\frac{\abs{M}}{r}\bigr)^{d-2}r^{d-2}\! ,
$\\[2mm]
where $\vrapprox{w}{\pm\eps}{a'}{b'}$ and $\vzapprox{w}{\pm\eps}{a'}{b'}$ are defined as in \eqref{volrdef} and \eqref{volzdef} respectively.
\end{lemma}
\bw
Using approximations in virtue of functions $g$ defined in \eqref{defg} we obtain by triangle inequality (Recall the definition of $\psi_{r,\pm\eps}$ in \eqref{defpsi_r}.)
\beqa
\Delta_{\pm\eps} &= &\Bigl|\int_{\Rd} I_{H_{r,M}}(x)\psi_{r,\pm\eps}(x)dx - \sum_{x\in\Zd} I_{H_{r,M}}(x)\psi_{r,\pm\eps}(x)\Bigr|\hskip3cm\nn\\
&\leq& \Bigl|\int_{\Rd}\bigl(I_{[a,b]} - g_{a',b',
w}\bigr)\bigl(Q[x-M]\bigr)\psi_{r,\pm\eps}(x)dx\Bigr|\nn \\ && \quad\quad\quad +\;\;
\Bigr|\vrapprox{w}{\pm\eps}{a'}{b'} -
\vzapprox{w}{\pm\eps}{a'}{b'}\Bigr|\nn\\ && \quad\quad\quad\quad\quad\quad+\;\;
\Bigl|\sum\limits_{x \in \Zd}\bigl(I_{[a,b]} - g_{a',b',
w}\bigr)\bigl(Q[x-M]\bigr)\psi_{r,\pm\eps}(x)\Bigr|. \eeqa Choosing
$a' , b' $ according to Lemma \ref{approx2} (1) and estimating the
first summand by taking the supremum, we obtain \beqa
\Delta_{\pm\eps}&\leq& \sup\limits_{{a' \in [a-w;a+w]}\atop{b' \in
[b-w;b+w]}}\!\Bigl|\int_{\Rd}\bigl(I_{[a,b]} - g_{a',b',
w}\bigr)\bigl(Q[x-M]\bigr)\psi_{r,\pm\eps}(x)dx\Bigr| \nn\hskip2cm\\
&& \hskip1cm +  \sup\limits_{{a' \in [a-w;a+w]}\atop{b' \in
[b-w;b+w]}}\!\Bigr|\vrapprox{w}{\pm\eps}{a'}{b'} -
\vzapprox{w}{\pm\eps}{a'}{b'}\Bigr|. \eeqa The application of Lemma
\ref{approx2} (2) completes the proof. \bwend
\vskip-10mm \phantom{....}


Repeating the proofs of Lemma \ref{approx1}, \ref{approx2} and \ref{approx3} in the case $a=a_r$ using Lemma \ref{volestdistr} instead of Lemma \ref{vol3} and \ref{vol4} we get immediately

\begin{lemma}\label{estmaindistr}
For $ F_{r,M}(b)$ defined as in \eqref{defFrMb} there exist constants $r_0=r_0(Q,M,b)\geq 1$ and $c_{Q,1}, c_{Q,2}\geq 1 $ depending on $Q$ and $d$ only, such that for $w>0,0<\eps<\frac{1}{4}$ the following estimate holds:
\beqa
\bigl|\volu_{\mathbb{Z}}F_{r,M}(b)-\volu F_{r,M}(b)\bigr|\leq \bigl(c_{Q,1}(b-a_r)\eps + c_{Q,2}w(1 +\eps)^{d-2}\bigr)r^{d-2}\phantom{xxxx} \nn\\\quad +\, \sup\limits_{{a' \in [a_r-w;a_r+w]}\atop{b' \in [b-w;b+w]}}\!\Bigr|\vrapprox{w}{\pm\eps}{a'}{b'} - \vzapprox{w}{\pm\eps}{a'}{b'}\Bigr|, \nn
\eeqa
where $\vrapprox{w}{\pm\eps}{a'}{b'}$ and $\vzapprox{w}{\pm\eps}{a'}{b'}$ are defined as in \eqref{volrdef} and \eqref{volzdef} respectively.
\end{lemma}
\begin{lemma}\label{pospartintegral}
For any $\beta >0, T \in \R$ it holds
\beqa
\ffrac{1}{2 \4\pi\4 i} \int^{\beta+ i\, \infty}_{\beta-i\, \infty}
\exp\{z\4T\}\4 \ffrac {dz} {z^2}= \max\{\4T,\40\4\} = T_+ .
\eeqa
\end{lemma}
\bw Complement the interval
$(\beta -i \, \infty, \beta + i \, \infty)$
by an infinite half circle in $\Rea z \ge0$
(resp. $\Rea z \le 0$) for $ T < 0$ (resp. $ T \ge 0$) and apply standard residue calculus.
\bwend

\begin{lemma}\label{thetasums1}
For a symmetric, $d\times d$ complex matrix $\Omega$, whose
imaginary part is positive definite the following holds:
\beqa\sum_{m \in \Z^d}\exp\Bigl[\pi i \cdot\Omega [m]+2\pi i\langle
m,v\rangle\Bigr]
&=&\Bigl(\det\Bigl(\frac{\Omega}{i}\Bigr)\Bigr)^{-\frac{1}{2}}\cdot
\exp\Bigl[ -\pi i \cdot\Omega^{-1}[ v]\Bigr]\nonumber\\ &\times&
\sum_{n \in \Z^d}\exp\Bigl[-\pi i \cdot\Omega^{-1} [n]+2\pi i\langle
n,\Omega^{-1}v\rangle\Bigr] \nn\eeqa and \beqa\int_{\R^d}\exp\bigl[\pi
i \cdot\Omega [x]+2\pi i\langle x,v\rangle\bigr]dx\nn\hskip3mm
&=&\left(\det\left(\frac{\Omega}{i}\right)\right)^{-\frac{1}{2}}\cdot
\exp\left[ -\pi i \cdot\Omega^{-1}[ v]\right], \eeqa where
$\Omega^{-1}\bigl[ x\bigr]$ denotes the quadratic form
$\langle\Omega^{-1}x,x\rangle$, defined by the inverse operator 
$\Omega^{-1}:\C^d\rightarrow\C^d$ (which exists since $\Omega$ is an
element of Siegel's upper half plane).\end{lemma} \bw See
\cite{mumford:1983}, p. 195 (5.6) and Lemma 5.8. \bwend
\begin{corollary}\label{thetasums1a} For $ z\in \C^d, \Rea z >0, v \in \C^d$ and a positive definite, symmetric $d\times d$ matrix $\Omega$ it holds
$$
\sum_{m \in \Z^d}\exp\bigl[ -z\Omega [m] +2\pi i\langle m,v\rangle\bigr] = \left(\det\left(z\cdot\frac{\Omega}{\pi}\right)\right)^{-\frac{1}{2}}\cdot\sum_{n\in \Z^d}\exp\left[-\frac{\pi^2}{z}\Omega^{-1}[n+v]\right].
$$
\end{corollary}
\bw Apply Lemma \ref{thetasums1} to the matrix $\frac{i}{\pi}z\Omega$.\bwend

\begin{lemma}\label{thetasums2}
For $z = \frac{1}{r^2} +i t, r>0, t \in \R$ and all $v \in \C^d$ it holds\\[2mm]

$ \sum\limits_{m \in \Z^d}\exp\Bigl[-\frac{2}{r^2}Q_+[m]-zQ[m] +2\pi i\langle m,v\rangle\Bigr]$\beqa & =& \det\Bigl(\frac{1}{\pi}\bigl(\frac{2}{r^2}Q_++zQ\bigr)\Bigr)^{-\frac{1}{2}}\cdot \exp\Bigl[ -\pi^2\bigl(\frac{2}{r^2}Q_++zQ\bigr)^{-1}\bigl[ v\bigr]\Bigr]\nonumber\\[2mm]
&&\quad\quad\times \sum_{n \in \Z^d} \exp\Bigl[ -\pi^2\bigl(\frac{2}{r^2}Q_++zQ\bigr)^{-1}\bigl[ n\bigr]  -2\pi^2\langle\bigl(\frac{2}{r^2}Q_++zQ\bigr)^{-1}n,v\rangle\Bigr],\nn\\[-2mm]&&\eeqa
and\\

$
\int_{\R^d}\exp\Bigl[-\frac{2}{r^2}Q_+[x]-zQ[x] +2\pi i\langle
x,v\rangle\Bigr]dx$
\beqa &=&\det\Bigl(\frac{1}{\pi}\bigl(\frac{2}{r^2}Q_++zQ\bigr)\Bigr)^{-\frac{1}{2}}\cdot
\exp\Bigl[ -\pi^2\bigl(\frac{2}{r^2}Q_++zQ\bigr)^{-1}\bigl[
v\bigr]\Bigr],\label{thetaintegral} \eeqa where
$\bigl(\frac{2}{r^2}Q_++zQ\bigr)^{-1}\bigl[ x\bigr]$ denotes the
quadratic form
$\langle\bigl(\frac{2}{r^2}Q_++zQ\bigr)^{-1}x,x\rangle$, defined by
means of the positive definite operator
$\bigl(\frac{2}{r^2}Q_++zQ\bigr)^{-1}:\R^d\rightarrow\R^d$.
\end{lemma}
\bw For $\Omega\defi \frac{i}{\pi} \left(\frac{2}{r^2}Q_+ + zQ\right)$ and $ z = \frac{1}{r^2} +i t, t \in \R$ the imaginary part $\Ima \Omega$ is positive definite. The application of Lemma \ref{thetasums1} to $\Omega$ completes the proof.\bwend
\begin{lemma}\label{thetaestimate1}
Let $\theta_v(z)$ and $\theta_{0,v}(z)$  denote the theta sum and
theta integral in \eqref{defthetasum} and \eqref{defthetaint}
respectively. Then there is a constant $c = c(Q)>0$, such that for \,$r \geq r_0 = r_0(Q,M)\geq 1$ and  $t
\in \R $, $|t|<\imbound$, \,the following bound holds \beqa \bigl|
\bigl(\theta_v - \theta_{0,v}\bigr)(r^{-2} + i\4 t) \bigr| & \ll_d &
q_0^{-\frac{d}{2}}\abs{r^{-2}+it}^{-\frac{d}{2}}\exp\left[-c\cdot
\Rea \bigl((r^{-2}+it)^{-1}\bigr)\right] \nn\eeqa $\hfill+ \;2 I_{(r,\infty)}(\abs{v}).$
\end{lemma}
\bw Using Lemma \ref{thetasums2} we obtain by \eqref{defthetasum}, \eqref{defthetaint} and the self-adjointness of the matrix $\bigl(\frac{2}{r^2}Q_+zQ\bigr)^{-1}$, that\\

$\bigl(\theta_v - \theta_{0,v}\bigr)(z)$\vskip-3mm\beqa &=& \exp\bigl[-zQ[M]\bigr]\det\Bigl(\frac{1}{\pi}\Omega\Bigr)^{-\frac{1}{2}}\cdot \exp\Bigl[ -\pi^2\Omega^{-1}\bigl[ -\frac{\tilde{v}}{2\pi r}\bigr]\Bigr]\nonumber\\
&&\quad\quad\quad\quad \quad\quad \quad\times \sum_{n \in \Z^d\setminus\{ 0\}} \exp\Bigl[ -\pi^2\Omega^{-1}\bigl[ n\bigr]  -2\pi^2\langle\Omega^{-1}n, -\frac{\tilde{v}}{2\pi r}\rangle\Bigr]\nn\\ &=&\exp\bigl[-zQ[M]\bigr]
\det\Bigl(\frac{1}{\pi}\Omega\Bigr)^{-\frac{1}{2}}\!\!\cdot\!\!\sum_{n \in \Z^d\setminus\{ 0\}}  \exp\Bigl[ -\Omega^{-1}\bigl[ \pi n-\frac{\tilde{v}}{2r}\bigr]\Bigr],\nn\\\label{thetaestimate1-eq1}
\eeqa
where $\Omega\defi \bigl(\frac{2}{r^2}Q_++zQ\bigr)$ and $\tilde{v}\defi 2r\Ima (z)QM - v$. Note that for $z = r^{-2} +it$ and $ t \leq \frac{1}{r} $ there exists a constant $c_0 = c_0(Q,M)>0$ such that $\abs{2r\Ima (z)QM} \leq c_0 $ uniformly in $r$. \\
Using $\det\left(\frac{1}{\pi}\Omega\right) = \frac{1}{\pi^d}\prod\limits_{1\leq j\leq d}\left(\frac{2}{r^2}|q_j|+zq_j\right)$ and $\Bigl|\frac{2}{r^2}|q_j|+zq_j\Bigr|  \geq \abs{q_j}\cdot\bigr| r^{-2}+it\bigl|$ for $z=r^{-2}+it$ and all $1\leq j\leq d$, we have
\beqa \left|\det\left(\frac{1}{\pi}\Omega\right)^{-\frac{1}{2}}\right| \leq \pi^{\frac{d}{2}}\cdot q_0^{-\frac{d}{2}}\cdot\abs{z}^{-\frac{d}{2}}.\label{thetaestimate1-eq2}
\eeqa

Since $\Omega$ can be orthogonal diagonalized, the matrix $  \Rea\left(\Omega^{-1}\right)$ has eigenvalues $\Rea \left(\bigl(\frac{2}{r^2}|q_j|+zq_j\bigr)^{-1}\right), 1\leq j\leq d$. For $t\leq\imbound$ we have
\beqa
\Rea \left(\bigl(\frac{2}{r^2}|q_j|+zq_j\bigr)^{-1}\right) \geq \frac{1}{|q_j|} \Rea (z^{-1})\geq \frac{1}{q} \Rea (z^{-1}),\;\; 1\leq j\leq d.
\nn\eeqa
Hence,
\beqa
\left|\exp\Bigl[ -\pi^2\Omega^{-1}\bigl[ n-\frac{\tilde{v}}{2\pi r}\bigr]\Bigr]\right| &=& \exp\Bigl[ -\Rea\left(\Omega^{-1}\bigl[\pi n-\frac{\tilde{v}}{2r}\bigr]\right)\Bigr]\nonumber\\ &=&
\exp\Bigl[ -\Rea\left(\Omega^{-1}\right)\bigl[\pi n-\frac{\tilde{v}}{2r}\bigr]\Bigr]\nonumber\\
&\leq& \exp\Bigl[ -\frac{1}{q}\Rea (z^{-1})\cdot \left|\pi n-\frac{\tilde{v}}{2r}\right|^2\,\Bigr].\label{thetaestimate1-eq3}
\eeqa
Using \eqref{thetaestimate1-eq1}, \eqref{thetaestimate1-eq2} and \eqref{thetaestimate1-eq3} we get\\

$\bigl| \bigl(\theta_v - \theta_{0,v}\bigr)(r^{-2} + i\4 t) \bigr|\quad \ll_d \quad \exp\bigl[-\frac{1}{r^2}Q[M]\bigr]\;q_0^{-\frac{d}{2}}|r^{-2}+it|^{-\frac{d}{2}}$\beqa\phantom{...}\hskip1cm\times\quad\sum_{n \in \Z^d\setminus\{ 0\}}\!\!\exp\Bigl[ -\frac{1}{q}\Rea \bigl((r^{-2}+it)^{-1}\bigr)\cdot \left|\pi n-\frac{\tilde{v}}{2r}\right|^2\,\Bigr].\label{expest2413}
\eeqa
For $ \abs{\tilde{v}}\leq \pi r$ we obtain \beqa\exp\Bigl[ -\frac{1}{q}\Rea \bigl((r^{-2}+it)^{-1}\bigr)\cdot \left|\pi n-\frac{\tilde{v}}{2r}\right|^2\,\Bigr]\leq \exp\Bigl[ -\frac{1}{q}\Rea \bigl((r^{-2}+it)^{-1}\bigr)\cdot \frac{\left|\pi n\right|^2}{2}\,\Bigr]\nn\eeqa and hence, for an appropriate constant $c = c(Q)>0$\\

{\large $
\sum\limits_{\text{\scriptsize $n \in \Z^d\!\setminus\!\{ 0\}$}}$}$\exp\Bigl[ -\ffrac{1}{q}\Rea \bigl((r^{-2}+it)^{-1}\bigr)\cdot \left|\pi n-\ffrac{\tilde{v}}{2r}\right|^2\,\Bigr]$\beqa &\leq&\sum_{n \in \Z^d\setminus\{ 0\}}\exp\Bigl[ -\frac{1}{q}\Rea \bigl((r^{-2}+it)^{-1}\bigr)\cdot \frac{\left|\pi n\right|^2}{2}\,\Bigr]\hskip3cm\nn\\ &\ll &\exp\Bigl[ - c\cdot\Rea \bigl((r^{-2}+it)^{-1}\bigr)\Bigr].
\label{thetaestimate1-eq4}
\eeqa
For $\abs{\tilde{v}}>\pi r$ set $\tilde{v} = L\pi r + w$, with $ L\in \Z , \abs{w}\leq \pi r$, then $w = \widetilde{v'}$ for $v' \defi v+ L\pi r$. By \eqref{thetaestimate1-eq1} we have obviously $\theta_v = \theta_{v'}$ and therefore we get by \eqref{thetaestimate1-eq4} and \eqref{expest2413} the inequality\\

$\bigl| \bigl(\theta_v - \theta_{0,v}\bigr)(r^{-2} + i\4 t) \bigr|$\beqa & \leq&\bigl| \bigl(\theta_{v'} - \theta_{0,v'}\bigr)(r^{-2} + i\4 t)\bigr| +\bigl| \bigl(\theta_{0,v'} - \theta_{0,v}\bigr)(r^{-2} + i\4 t)\bigr|\nn\\ &\ll_d& \exp\bigl[-\frac{1}{r^2}Q[M]\bigr] \; q_0^{-\frac{d}{2}}|r^{-2}+it|^{-\frac{d}{2}}\cdot \exp\Bigl[ - c\cdot\Rea \bigl((r^{-2}+it)^{-1}\bigr)\Bigr] \nonumber\\&&\hskip50mm +\;\; \bigl|  \theta_{0,v'}(r^{-2} + i\4 t)\bigr|+ \bigl|  \theta_{0,v}(r^{-2} + i\4 t)\bigr|\nn\\&\leq& \exp\bigl[-\frac{1}{r^2}Q[M]\bigr]\; q_0^{-\frac{d}{2}}|r^{-2}+it|^{-\frac{d}{2}}\cdot\exp\Bigl[ - c\cdot\Rea \bigl((r^{-2}+it)^{-1}\bigr)\Bigr] + 2.\nn\\
\label{thetaestimate1-eq5}
\eeqa 
The result now follows by \eqref{thetaestimate1-eq4}, \eqref{expest2413} and \eqref{thetaestimate1-eq5} for $r \geq r_0$, $ r_0\geq 1$ large enough, since $\abs{\tilde{v}}>\pi r $ implies $ \abs{v}\geq \pi r - c_0(Q,M)\geq r$ for $r$ large enough.
\bwend
\begin{lemma}\label{thetaestimate2}
Let \,$\theta_v(z)$ denote the theta function in 
\eqref{defthetasum} depending on $\Q$ and $v\in \C^d$.
\,For \,$r \geq r_0 =r_0(Q,M)\geq 1 $, \,$t \in {\R}$, \,the following bound holds
\beqa
\bigl| \theta_v (r^{-2} + i\4 t) \bigr| & \ll_d & (\det \Omega)^{-1/4} \4
r^{d/2}\,
\psi (r,t)^{1/2},\; \text{ where }\label{eq:thetaestimate2}\\
\psi (r,t) & = & \sum_{m,n \in \mathbb{Z}^d} \exp \Bigl\{- \frac{r^2}{2}\4
 \Omega^{-1} [\4\pi\4m - 2\4t\, \Q  n\4] - \frac{2}{r^2}\4 \Omega[\4n\4]\Bigr\},
\nn\eeqa with $\Omega \defi 2\!\cdot\!\Q_+ +\Q$.
\end{lemma}

Note that the right hand side of this inequality is independent of $v \in
{\C}^d$.\\

\bw
For any $x, y \in {\R}^d$ the equalities
\beqa
2\, \left(\Omega[\4x\4] + \Omega[\4y\4]\right) & = & \Omega[\4x + y\4] + \Omega[\4x -y\4], \\[0.2cm]
\left\langle \Omega\4 (x + y), \,x - y \right\rangle
& = & \Omega[\4x\4] - \Q[\4y\4] \label{eq:quad1}
\eeqa
hold.
Rearranging $\theta_v(z)\, \overline{\theta_v(z)}$
and using \eqref{eq:quad1}, we would
like
 to use $m + n$ and $m-n$ as new summation variables on a
lattice. But both vectors have the same parity, i.e., $m + n \equiv
m - n \mod 2$. Since they are dependent one has to consider the
$2^d$ sublattices indexed by \,$ \alpha = (\alpha_1, \dots, \alpha_d)$
\,with
 \,$\alpha_j = 0,1$, \,for \,$1 \leq j \leq d$:
\beqa
\mathbb{Z}^d_{\alpha}& \defi & \{ m \in \mathbb{Z}^d \,:\, m \equiv \alpha
\mod 2 \},\nn
\eeqa
where, for $ m = (m_1, \dots , m_d)$,
\,$m \equiv \alpha \mod 2$ \,means
 \,$m_j \equiv \alpha_j \mod 2$, \,$ 1\leq j \leq d$.
\,Thus writing
$$
\theta_{v,\alpha}(z) \defi \sum_{m \in \mathbb{Z}^d_{\alpha}} \exp\left[-\frac{2}{r^2}Q_+[m] - z\cdot Q[m] - i\cdot\langle m,\frac{v}{r}-2tQM\rangle\right],
$$
we obtain $\theta_v(z) =  \exp\bigl[-zQ[M]\bigr]\sum_{\alpha} \theta_{v,\alpha}(z)$ and hence by
the Cauchy--Schwarz inequality
\beqa
\bigl| \theta_v(z) \bigr|^2 \leq \ 2^d  \exp\bigl[-\frac{2}{r^2}Q[M]\bigr]\,\sum_{\alpha}\
\bigl|\theta_{v,\alpha}(z) \bigr|^2.\label{eq:cauchy}
\eeqa
 Using \eqref{eq:quad1}
and the absolute convergence of \,$\theta_{\alpha}(z)$, \,we
may rewrite the quantity $\theta_{v,\alpha}(z)\, \overline{\theta_{v,\alpha}(z)}$ for $ z = \frac{1}{r^2}+it $ and $ \tilde{v} \defi v-2trQM $ as\\

$\theta_{v,\alpha}(z) \, \overline{\theta_{v,\alpha}(z)}$\beqa
 &=& \sum_{m, n \in \Z^d_{\alpha}} \exp\left[-\frac{1}{r^2}\bigl(\Omega[m]+\Omega[n]\bigr)-it\cdot\bigl(\Q[m]-\Q[n]\bigr) - i\cdot \langle m-n,\frac{\tilde{v}}{r}\rangle\right]\nn\\
&=& \sum_{m, n \in \Z^d_{\alpha}} \exp\left[-\frac{2}{r^2}\bigl(\Omega[\ov{m}]+\Omega[\ov{n}]\bigr)-2i\cdot \bigl\langle 2t\cdot\Q\ov{m}+\frac{\tilde{v}}{r},\ov{n}\bigr\rangle\right]
\eeqa
where \,$\ov{m} = \ffrac{m+n}{2}$, \,$ \ov{n} = \ffrac{m - n}{2}$. \\ Note that the map
$
 H : \bigcup_{\alpha}\mathbb{Z}^d_{\alpha} \times \mathbb{Z}^d_{\alpha} \rightarrow
\mathbb{Z}^d \times \mathbb{Z}^d,
(m, n)  \longmapsto \Big(\ffrac{m + {n}}{2}, \ffrac{m-n}{2}\Big)
$ is a bijection. Therefore we get by \eqref{eq:cauchy}\\

$\exp\bigl[\frac{2}{r^2}Q[M]\bigr]\cdot\bigl| \theta_v(z) \bigr|^2 $\beqa &\ll_d& \sum_{\alpha \in\{0,1\}^d}\sum_{\ov{m}, \ov{n} \in \Z^d_{\alpha}} \exp\left[-\frac{2}{r^2}\bigl(\Omega[\ov{m}]+\Omega[\ov{n}]\bigr)-2i\cdot \bigl\langle 2t\cdot\Q\ov{m}+\frac{\tilde{v}}{r},\ov{n}\bigr\rangle\right]\nn\\ &=& \sum_{\ov{m}, \ov{n} \in \Z^d} \exp\left[-\frac{2}{r^2}\bigl(\Omega[\ov{m}]+\Omega[\ov{n}]\bigr)-2i\cdot \bigl\langle 2t\cdot\Q\ov{m}+\frac{\tilde{v}}{r},\ov{n}\bigr\rangle\right].\label{eq:double1}
\eeqa
In this double sum fix $\ov{n}$ and sum over $\ov{m}
\in \mathbb{Z}^d$ first. Using Corollary \ref{thetasums1a} for $ z = \frac{2}{r^2},$ we get for $\delta\defi \left( \det \left(\frac{2}{\pi r^2}\cdot \Omega\right)\right)^{-\frac{1}{2}}$ by the symmetry of $\Q$
\beqa
\theta_v(z, \ov{n}) & \defi & \sum_{\ov{m} \in \Z^d }\exp\left[-\frac{2}{r^2}\bigl(\Omega[\ov{m}]+\Omega[\ov{n}]\bigr)-2i\cdot \bigl\langle 2t\cdot\Q\ov{m}+\frac{\tilde{v}}{r},\ov{n}\bigr\rangle\right]\nn \\
& = &\hskip-3mm
 \delta
\sum_{m \in \mathbb{Z}^d}
 \exp\left[ - \frac{\4 r^2}{2}\4
\Omega^{-1}[\4 \pi\4 m - 2\4t\,\Q \4\ov{n}\4]-\frac{2}{r^2}\Omega[\ov{n}]-2i\langle\4\frac{\tilde{v}}{r},\ov{n}\4\rangle\right].\nn
\eeqa
Thus,
\beqa
 \bigl| \theta_v(z, \ov{n}) \bigr| \leq
 \delta
\4 \sum_{m \in \mathbb{Z}^d}
 \exp \Bigl\{ - \frac{r^2}{2}\4
 \Omega^{-1}[\4\pi\4 m- 2\4t\,\Q\4 \ov{n}\4]-\frac{2}{r^2}\Omega[\ov n ]\Bigr\}.\quad \label{eq:double2}
\eeqa
Hence, we obtain by \eqref{eq:double1} and \eqref{eq:double2}\\[2mm]
$
\bigl| \theta_v(z)\bigr|^2 \ll_d  \exp\bigl[-\frac{2}{r^2}Q[M]\bigr]\;\;(\det \Omega)^{-1/2}\4 r^d $\\
$$ \hskip3cm \times\;\sum_{m,\ov n \in
 \mathbb{Z}^d} \4 \exp \Bigl\{- \frac{r^2}{2}\4 \Omega^{-1}[\4\pi\4m
- 2\4 t\, \Q\4\ov{n}\4] - \frac{2}{r^2}\4 \Omega[\4\ov{n}\4]\Bigr\},
$$
which proves Lemma \ref{thetaestimate2}  for $r> r_0 = r_0(Q,M)\defi\abs{Q[M]}^{1/2}+ 1$   .
\bwend

In the following we shall use some facts in the geometry of numbers
(see \cite{davenport:1958}).\\[2mm] Let $F: {\R}^d \rightarrow [\40, \infty\4)$
denote a norm on ${\R}^d$, that is $F( \alpha\4 x) =
\lvert \alpha \rvert\, F(x),$ for $\alpha \in {\R}$,
 and $F(x + y) \leq F(x) +F(y)$.
The successive minima $M_1 \leq \dots \leq M_d$ of $F$ with
respect to the lattice $\mathbb{Z}^d$ are defined as follows: Let
\,$M_1 = \inf\left \{F(m): m \neq 0, m \in \mathbb{Z}^d \right\}$ \,and define
$M_k$ as the infimum of $\lambda > 0$
such that the set \,$\left\{m \in \mathbb{Z}^d \,:\, F(m) < \lambda\right \}$
\,contains $k$ linearly independent vectors. It is
easy to see that these infima are attained, that is there exist
linearly independent vectors $a_1, \dots, a_d \in \mathbb{Z}^d$ such
that $F(a_j) = M_j$.

\begin{lemma} \label{Davenport}
Let \,$ L_j (x) = \sum^d_{k = 1} q_{jk}\4 x_k$, $1 \leq j \leq d$, \,
denote linear
forms on ${\R}^d$ such that $q_{jk} = q_{kj}$, \,$j,k=1,\dots ,d$.
Assume that $r \geq 1$ and let $\lVert v \rVert $
denote the distance of the number
$v$ to the nearest integer. Then the number of $m=(m_1, \dots,
m_d) \in \mathbb{Z}^d$ such that
\beqa
\norm{L_j(m)} < r^{-1}, \; \lvert m_j \rvert < r,
\qquad\mbox{for all} \ \ 1 \leq j \leq d,\nn
\eeqa
 is bounded from above by \,$c_d\4 (M_1
 \cdots M_d)^{-1} $, \,
where $c_d>0$ denotes a constant depending on\/~ $d$ only,
$M_1 \leq \dots \leq M_d$
 are the first $d$
 of the $2\4d$ successive minima $M_1 \leq \dots \leq M_{2d}$ of the
 norm $F: {\R}^{2d} \rightarrow [\4 0, \infty)$ defined for
 vectors \,$y = (x, \ov{x}) \in {\R}^{2d}$, \,$ x, \ov{x}
 \in {\R}^d$, \, $\ov{x} =(\ov{x}_1, \dots, \ov{x}_d)$, \, as
\beqa
F(y) \defi \max \bigl\{ r\4 \lvert\4 L_1(x)- \ov{x}_1\4\rvert ,
\dots,r\4 \lvert\4 L_d(x)- \ov{x}_d \4\rvert,
\,r^{-1}\4 \lvert x \rvert_\infty\bigr \}.
\label{eq:norm}
\eeqa
Moreover,
\beqa
\frac{1}{2\4d} \leq M_k\4 M_{2d+1-k} \leq(2\4d)^{2d-1}, \qquad
1 \leq k \leq 2\4d.\label{eq:Dbound}
\eeqa
\end{lemma}

\bw  \cite{davenport:1958}, (20), p. 113, Lemma 3.\bwend

Note that for some constant, say $c(d) >0$,  depending on $d$ only
\beqa \label{eq:minima0}
 r^{-1} \le M_1\le \dots\le M_d\le c(d) ,
\eeqa
where the first inequality is obvious by
\,$F(m, \ov m) \geq r^{-1}\4 \abs{ m }_{\infty}$.
\, If here\/ \,$m=0$ \,then \/$\ov m \ne 0$
\, and \, $F(m,\ov m )= r\4\lvert \ov m \rvert_\infty\ge
r^{-1}\4\lvert \ov m \rvert_\infty \ge r^{-1}.$
Finally,  \,$M_d \ll_d 1$ \,follows from
\eqref{eq:Dbound}
 for \,$k = d$.\\

In the following we shall consider linear forms
\beqa
L_j(x) = \sum^d_{k=1} t\, q_{jk}\,x_k, \qquad 1 \leq j \leq d,
\label{eq:forms}
\eeqa
where $\Q = (q_{ij})$, $i,j = 1, \dots, d$, denotes the components of the
matrix $\Q$ and where $t \in {\R}$ is
arbitrary. We denote the corresponding successive minima of the norm
$F(\cdot )$ defined by \eqref{eq:norm} and
\eqref{eq:forms} for fixed $t$ by $M_{j,t}$, \,$ j = 1,\dots, d$.
\,Thus, we can write
\beqa\label{eq:Mjt}
M_{j,t}=\bigl|L(m,n,t)\bigr|_\infty,
\eeqa
for some $m,n\in\Zd$, where
\beqa
L(m, n, t) = \bigl(r\4 (m_1 - t \4 (\Q\4 n)_1), \dots,
 r \4 (m_d - t\4(\Q\4 n)_d ),
\,r^{-1}\4 n_1, \dots, r^{-1}\4n_d\bigr). \nn\\[-1mm]\label{eq:defL}
\eeqa
It is easy to see from the definition that
\beqa
M_{j,t}=M_{j,-t}, \qquad j = 1,\dots, d ,\; t\in\R.
\label{eq:symm}
\eeqa
\vskip 0.2cm

\begin{lemma}\label{thetaestimate3} Let \,$r \geq 1$. \,Then
\beqa
\bigl\lvert \theta (r^{-2} + i\,t\frac{\pi}{2})\bigr\rvert
\ll_d 
q_0^{-\frac{3d}{4}} r^{d/2}\4 (M_{1,t} \cdots M_{d,t})^{-1/2}.\nn
\eeqa
\end{lemma}

\bw
By Lemma \ref{thetaestimate2} we need to estimate
the theta series $\psi (r,t\4\pi/2)$. Since the matrix $ \Omega = 2\Q_++ Q$ is positive definite we may use the inequalities
\,$\Omega^{-1}[\4x\4] \ge \frac{1}{3q}\4 \abs{ x }^2_{\infty}$
\,and \,$\Omega[\4x\4] \ge q_0\abs{x}^2_{\infty}$,  and we get
with  \,$c_\Q=\min\left\{\frac{\pi^2}{6\,q},2q_0\right\}$
\beqa
\psi (r, t \,\frac{\pi}{2}) \ll_d \sum_{m, n \in \mathbb{Z}^d} \exp \bigl\{-
c_\Q\,\lvert L(m,n,t)\rvert^2_{\infty} \bigr\},\, \label{eq:psi0}
\eeqa
where $L(m, n, t)$ is defined in \eqref{eq:defL}. Let
$$H\defi \bigl\{(m, n) \in \mathbb{Z}^{2d}\, : \,
\lvert L(m, n, t) \rvert_{\infty} < 1 \bigr\}. $$ Now,
Lemma \ref{Davenport} may be restated for the forms \eqref{eq:forms} as
\beqa
\# H \ll_d (M_{1,t} \cdots M_{d,t})^{-1}.
\eeqa

In order to bound \,$\psi (r,t\4\pi/2)$, we introduce for
\,$k \defi (k_1, \dots, k_{2d})\in \Z^{2d}$ \,the sets
\beqa
B_k & \defi & \left[k_1- \ffrac{1}{2}, k_1
+ \ffrac{1}{2}\right)\times\cdots\times \left[k_{2 d}-
\ffrac{1}{2}, k_{2 d} + \ffrac{1}{2}\right) \; \text{ and } \nonumber \\[0.2cm]
H_k & \defi &
\bigl\{ (m, n) \in \mathbb{Z}^{2d}: L(m,n,t) \in B_k \bigr\}\nonumber
\eeqa
such that ${\R}^{2d} = \bigcup_{k} B_k$.
For any fixed \,$(m^*,n^*) \in H_k$ \,we have
\beqa
(m - m^*, n - n^*)\in H\qquad\text{for any}\;(m,n)\in H_k.
\nonumber
\eeqa
Hence, we conclude for any $k\in \Z^{2d}$
\beqa
\# H_k \, \leq \, \# H \, \ll_d \, (M_{1, t} \cdots
M_{d,t})^{-1}.\label{eq:box1}
\eeqa

Since \,$x \in B_k $ \,implies \,
$ \abs{x}_{\infty} \geq \abs{ k }_{\infty}/2$, \,
we obtain by \eqref{eq:psi0} and \eqref{eq:box1}
\beqa
\psi (r,t\4\pi/2) & \ll_d & \#H_0 +
 \sum_{k \in \Z ^{2d} \setminus 0}
\; \sum_{m,n \in \Z ^{2d}} \I\bigl\{L(m,n,t) \in B_k\bigr\}\4
 \exp \bigl\{-c_\Q\4 \abs{ k }^2_{\infty}/4\bigr\} \nonumber \\
 & \ll_d & (M_{1,t} \cdots M_{d,t})^{-1} \sum_{k \in \Z ^{2d}}
 \4 \exp \bigl\{- c_\Q\4 \abs {k}_{\infty}^2/4\bigr\}
\nonumber \\
 & \ll_d &(M_{1,t} \cdots M_{d,t})^{-1} (c_\Q^{-1/2}+1)^{2d},\nn
\eeqa
using similar bounds as in \eqref{thetaestimate1-eq4}. Some simple bounds together
with Lemma \ref{thetaestimate2} finally conclude the proof of Lemma \ref{thetaestimate3}.
\bwend

In the following we consider an arbitrary, real,  symmetric, non-degenerate $d^*\times d^*$  - matrix $\Q^*$. The norm on $\R^{d^*}$\!, associated by \eqref{eq:defL}, and the associated successive minima are denoted by $L^*$ and $ M^*_{j,t},1\leq j\leq d^*$, respectively.

\begin{lemma}\label{gaps}
Let \,$(m,n), (m', n') \in \Z ^{2d^*} \setminus 0;$  \,\,$ t, t' \in \R$ and
\,$r\ge1$. \,Let
$M\defi \lvert L^*(m,n,t)\rvert_{\infty}$ and $M' \defi
\lvert L^*(m', n', t')\rvert_{\infty}$. Assume
that $\langle\Q^* n,n' \rangle  >0$ and
\beqa\label{eq:Mbound}
\max \{\4M, M'\4\} \le (4\,d^*)^{-1/2}.
\eeqa
Then for
\beqa\label{eq:Delta}
\Delta=\Delta(m,n;m',n')\defi
 \bigl\lvert \langle\4 n',m\4\rangle - \langle\4 m',n\4\rangle \bigr\rvert
\eeqa
 the following holds:

\beqa
{\rm i)}\;\;\;\Delta = 0 &\Rightarrow& \lvert\4 t-t'\rvert \leq \,
\ffrac{(d^*)^{1/2}\4\max \{\4M, M'\4\} \, (\lvert n\rvert
 + \vert n' \rvert)}{r\,
 \langle\4 \Q^*\4n, n' \4\rangle },\nonumber \\
&&\label{eq:alter}\\
{\rm ii)}\;\;\;\Delta \neq 0&\Rightarrow&
\lvert\4 t-t'\rvert  \geq \,\langle\4 \Q^* \4n,n' \4\rangle ^{-1}\!/2.\nonumber
\eeqa

\vskip2mm
In particular, assuming  $n=n'$ and \eqref{eq:Mbound} the alternative {\rm i)} in \eqref{eq:alter}  holds.

Furthermore, assuming $(m,n) \in \Z ^{2d^*}\! \setminus 0$ and
$M=\lvert L^*(m,n,t)\rvert_{\infty} \le (4\,d^*)^{-1/2} $
we have either
\beqa\label{eq:alter0}
\mbox{\rm i)} \quad \lvert t\rvert
\leq \ffrac{2\,d^*\4 M \4 \abs{n}}{r\, \abs{ \Q^*\4 n}}
 \quad or \quad
\mbox{\rm ii)} \quad \abs{t} \geq \ffrac{1}{2\,\abs{ \Q^*\4 n}}.
\eeqa
This means $t,t'$ resp. $t,0$
 have to be either 'near'
to each other or 'far' apart.
\end{lemma}

\bw \cite{goetze:2004}, p. 217, Lemma 3.6 or \cite{elsner:2006}, p. 38, Lemma 2.4.17 \bwend

The application of Lemma \ref{gaps} to $\Q^* = \Q$ yields the following

\begin{corollary}\label{multineq}
Let $r \geq 1$ and $d\ge 4$. Then
\beqa \label{eq:multbound}
 M_{1,t} \cdots M_{d,t} \ge
 d^{-d}\Bigl(\min \Bigl\{\ffrac{q_0\lvert t\rvert\4 r}2,\,
\ffrac 1 {q\4 \lvert t \rvert\4 r}\Bigr\}\Bigr)^d.
\eeqa
\end{corollary}

\bw Since $\lvert \Q\4 n\rvert =\lvert \Q_+\4 n\rvert $ we have $\lvert \Q\4 n\rvert \geq q_0\lvert n \rvert $,
and $ \lvert n\rvert \geq q^{-1}\vert \Q\4 n \rvert $. In the case, where $ M_{j,t}\le(4\,d)^{-1/2}$ we obtain by \eqref{eq:alter0}, $\lvert n \rvert_{\infty} \leq r\,M_{j,t}$ \, and
\,$2\4d^{1/2} \leq d $:
\beqa
\mbox{\rm i)}& &\lvert t \rvert \4 r\4 d^{-1}q_0 \leq \lvert t \rvert
\4 r\4 d^{-1}\4 \ffrac{\lvert \Q\4 n \rvert } {\vert n \rvert}
\leq 2\4M_{j,t} \nonumber \\[-0.1cm]
\text{ or \qquad} &&\label{eq:3.30} \\[-0.1cm]
\mbox{\rm ii)}
& &\ffrac{1}{\lvert t \rvert } \leq 2\,\lvert \Q\4 n\rvert \leq
2\4q\4 \vert n\rvert \leq
2\4d^{1/2}\4q\4 \lvert n \rvert_{\infty}
\leq q\4 d\4 r\4 M_{j,t},\nonumber
\eeqa
for appropriate $(m,n) \in \Z ^{2d}$ depending on $j$ such that
$M_{j,t} = \lvert L(m,n,t)\rvert _{\infty}$.
Note that if \,$ M_{j,t}\ge (4\,d)^{-1/2}$,
\,then \,$ M_{j,t}\ge d^{-1}$ \,since \,$d\ge 4$.
\,Combined with \eqref{eq:3.30},
this proves
 Corollary \ref{multineq} since
\beqa \label{eq:multbound1}
\min \Bigl\{\ffrac{q_0\lvert t\rvert\4 r}2,\,
\ffrac 1 {q\4 \lvert t \rvert\4 r}\Bigr\}\le1\nonumber
\eeqa
(recall that $q_0\leq q$).
\bwend

\vskip3mm In the following two Lemmas we will additionally assume
that the matrix ${\Q^*}$ is positive definite. The smallest and the
largest eigenvalue of ${\Q^*}$ is denoted by $q_0^*$ and $q^*$
respectively.

\begin{lemma}\label{integral}
Let $[\4\kappa,\xi\4] \subset \R, \,\, 0 < \kappa < \xi < \infty$.
Define for $g \in C^1 [\4\kappa ,\xi\4]$
such that
$g \geq 0, \,\, g' \leq 0$~ on ~$[\4\kappa, \xi\4]$,
\beqa
H_{\kappa,\xi}(\tau)\defi H_{\kappa,\xi, \Q^*}(\tau) \defi \int^\xi_\kappa \4\I\{M^*_{1,t} \leq \tau\}\, g (t)\, dt.
\label{eq:Hdef}
\eeqa
 Then, for all
\beqa\label{eq:taubound}
 \kappa > \left(q^*_0r\right)^{-1},\qquad
r^{-1} \leq \tau \leq (2\4 d^*)^{-1},
\eeqa
we have
\beqa
\label{eq:integral} H_{\kappa,\xi} (\tau) &\ll_{d^*} & \ovln H_{\kappa,\xi}(\tau)\defi
 \ffrac{q^*}{q^*_0} \4\tau^2
\int^\xi_{\kappa(\tau r)}g(t)\, dt +\ffrac{1}{q_0^*}\ffrac \tau r g(\kappa(\tau\4r)),
\eeqa
where \,$\kappa(v) = \max\bigl\{\kappa, (2\4q^*\4v\4d^{1/2})^{-1}\bigr\}$,
\, provided that \, $\kappa(\tau\4r) \le  \xi$.
\,In the case where \,$\kappa(\tau\4r) > \xi$,
 \,we have \,$ H_{\kappa,\xi} (\tau) =0$.
\end{lemma}

\bw \cite{goetze:2004}, p. 219, Lemma 3.8 \bwend

For indicator functions $g$  Lemma \ref{integral}
 reads as follows.

\begin{lemma}\label{indicator}
Let
$\lambda$ denote the Lebesgue measure.
 There exists a constant
~$c(d^*)$ depending on $d^*$ only such that for any $r\ge 1$, $\tau >0$ and
any interval $[\kappa,\xi]$ with $\xi>\kappa$ the following holds:
\beqa\nonumber
I(\tau)\defi \lambda\{ t \in [\kappa,\xi] \, : \, M^*_{1,t} \le \tau\}\, \le\, c(d^*)\left(\ffrac{q^*}{q^*_0}\4 \tau^2 \4(\xi-\kappa) + \ffrac{1}{q^*_0}\tau\4 r^{-1}\right).
\eeqa
\end{lemma}

\bw \cite{goetze:2004}, p. 222, Lemma 3.9 \bwend

We now return to general (not necessary positive definite) non-degenerate, symmetric, real $d\times d$ - matrix $\Q$, to the corresponding norm $ L$ (see \eqref{eq:defL}) and the associated successive minima $M_{j,t}$(see \eqref{eq:Mjt}).\\[2mm]
In the sequel we will assume, that $\Q$ is a {\em block-type} matrix, that is, that there exist positive definite matrices $ \Q^+ \in {\rm GL}\bigl(\R^{d^+}\bigr) , \Q^- \in {\rm GL}\bigl(\R^{d^-}\bigr),\, d^+ + d^- \geq 5$ with
\beqa
\Q \;\;= \;\;\left(\begin{array}{cc}\Q^+& 0\\0& -\; \Q^-\end{array}\right) .
\nn\eeqa

We denote the corresponding successive minima of the norm
$F^{\pm}(\cdot)$, defined by the analogon of \eqref{eq:norm} and
\eqref{eq:forms} for $\Q^{\pm}$, for a fixed $t$ by $M^{\pm}_{j,t}$,
\,$ j = 1,\dots, d^{\pm}$. \,Thus, we can write
\beqa\label{eq:Mjtpm}
M^{\pm}_{j,t}=\left|L^{\pm}(m,n,t)\right|_\infty, \eeqa for some
$m,n\in\Z^{d^{\pm}}$, where \beqa L^{\pm}(m, n, t)\! =\! \left(r\4 (m_1
- t \4 (\Q^{\pm}\4 n)_1), \dots,
 r \4 (m_{d^{\pm}} - t\4(\Q^{\pm}\4 n)_{d^{\pm}} ),
\frac{1}{r}\, n_1, \dots, \frac{1}{r}\, n_{d^{\pm}}\right).\nn
\eeqa
As in \eqref{eq:symm} we have
\beqa
M^{\pm}_{j,t}=M^{\pm}_{j,-t}, \qquad j = 1,\dots, d^{\pm} ,\; t\in\R.
\label{eq:symmpm}
\eeqa

In this special case there is a simple relation between the first successive minimum of $\Q$ and those of $\Q^+$ and $\Q^-$.
\begin{lemma}\label{Mzerlegung}
For $t \in \R$ holds
\beqa
M_{1,t} \geq \min\left\{ M^+_{1,t}, M^-_{1,t} \right\}. \label{eq:Mzerlegung}
\eeqa
In particular, for $\tau \in \R$,
\beqa
\I\{M_{1,t}\leq \tau\} &\leq &\I\{M^+_{1,t}\leq\tau\} +
\I\{M^-_{1,t}\leq \tau\}.\nn
\eeqa
\end{lemma}
 \bw
Choose $(m,n) = \left({\small \left({m_+ \atop m_-}\right),\left({n_+ \atop n_-}\right)}\right) \in \Zd\setminus 0$ such that $M_{1,t} = \abs{L(m,n,t)}_{\infty}$. It is easy to see, that $$M_{1,t} = \abs{L(m,n,t)}_{\infty} = \max\left\{\bigl|L^+(m_+,n_+,t)\bigr|_{\infty},\bigl|L^-(m_-,n_-,-t)\bigr|_{\infty}
\right\}.
$$ Since $(m,n)\neq 0$, it follows $(m_+,n_+)\neq 0$ or  $(m_-,n_-)\neq 0$ and hence by \eqref{eq:symmpm},
$$\bigl|L^+(m_+,n_+,t)\bigr|_{\infty} \geq M^+_{1,t} \;\;\;\; \text{ or }\;\;\;\; \bigl|L^-(m_-,n_-,-t)\bigr|_{\infty} \geq M^-_{1,-t} =   M^-_{1,t}.$$
This proves \eqref{eq:Mzerlegung}.
\bwend

\begin{corollary}\label{lemmadiophantine1}
Again, $\lambda$ denotes the Lebesgue measure. Then
 there exists a constant
~$c = c(d)>1$ depending on $d$ only, such that for any $r\ge 1$,
$\tau
>0$ and any interval $[\kappa,\xi]$ with $\xi>\kappa$ the following
holds: \beqa\nonumber I(\tau)\defi \lambda\{ t \in [\kappa,\xi] \, :
\, M_{1,t} \le \tau\}\, \le\, c\cdot\left(\ffrac{q}{q_0}\4 \tau^2
\4(\xi-\kappa) + \ffrac{1}{q_0}\tau\4 r^{-1}\right). \eeqa
\end{corollary}
\bw Using Lemma \ref{indicator} and Lemma \ref{Mzerlegung} we obtain
\beqa I(\tau) &\leq& \int \I\{M^+_{1,t}\leq\tau\} +
\I\{M^-_{1,t}\leq \tau\}\,\lambda (dt) \nn\\ &\leq& \bigl(c(d^+)+
c(d^-)\bigr)\left(\ffrac{q}{q_0}\4 \tau^2 \4(\xi-\kappa) +
\ffrac{1}{q_0}\tau\4 r^{-1}\right),\nn\eeqa where we have used, that
$q$ (resp. $q_0$) is larger (resp. smaller) than the corresponding
largest (resp. smallest) eigenvalue of $Q^+$ and $Q^-$. Taking $c
\defi \max\limits_{{d^+, d^-\in \N}\atop{d^++d^- = d}}\bigl(c(d^+)+
c(d^-)\bigr)$ completes the proof.\bwend

\begin{lemma}\label{final}
Let ${M}(t) = M_{1,t} \cdots
M_{d,t}$,       $\gamma =\gamma({\kappa,\xi}) = r^d\4 \4
\inf_{\kappa \leq t \leq \xi} {M}(t)$ and introduce
\beqa
D= \max\bigl\{(2\4d)^{-d}\4 r^d,\,\gamma\bigr\}
 \quad
\text{ and } \quad G(\kappa,\xi) = \int_{\kappa}^{\xi}\4 g(t)\, dt,
\nn\eeqa
 for $0< \kappa < \xi\le\infty$ and let ~$g(t)$ and ~$\kappa(v)$ be
as in Lemma \ref{integral}.
For  $\kappa > \xi$ we define $G(\kappa,\xi) =0$.
 Then
\beqa
I_{\kappa,\xi}&\defi &\int_{\kappa}^{\xi}\ffrac{g(t)}{ M(t)^{1/2}}\,dt\nonumber \\
&\ll_d & q_0^{-1}r^{d/2-2}\int_{\gamma}^{D} v^{-1/2+1/d}\left(q \4 v^{1/d}
G(\kappa(v^{1/d} ),\xi) \, +\,  g(\kappa(v^{1/d})) \right)\ffrac{dv}{v} \nn\\ &&+ \quad G(\kappa,\xi).\label{eq:finalint}
\eeqa

\end{lemma}

\bw We generalize the proof in \cite{goetze:2004}, p. 222, Lemma 3.10:

Write \,$\ov{\gamma}\defi  \inf_{\kappa \leq  t  \leq \xi} {M}(t)$ \, and
\,$c_d= (2\4d)^{-d}$. If $\ov\gamma\ge c_d$, then $I_{\kappa,\xi}\ll_d G(\kappa,\xi)$
and \eqref{eq:finalint} is obvious. In the case \beqa
\ov\gamma< c_d\label{eq:gammacd}\eeqa
 we define
\beqa\label{eq:Jdef}
J_{\kappa,\xi}(v)\defi \int^{\xi}_{\kappa} g(t) \,I_{\{{M}(t) \leq \, v\}}\, dt
\eeqa
for $0 < \kappa <\xi.$ Since $M_{j,t} \leq M_{d,t} \ll_d 1,$ for $j=1,\dots , d,$
 by Lemma \ref{Davenport}, there exists a constant  $\bar{M}$ depending on~$d$ only such that
${M}(t)\le \bar{M}$ for all ~$t$. Therefore we have for all $t\in [\4\kappa,\xi\4]$
$${M}(t)^{-1/2} = \int_{\ov{\gamma}}^{\bar{M}}
\eps^{-1/2} \4 d\4 I_{\{ {M}(t) \le \eps\}} .
$$
Hence, Fubini's Theorem implies
$$
I_{\kappa,\xi}= \int_{\ov{\gamma}}^{\bar{M}} \eps^{-1/2} \4 d J_{\kappa,\xi}(\eps).
$$
Splitting the integral $I_{\kappa,\xi}$ into the
part where $\eps \le c_d$ and its complement, we obtain
\beqa
I_{\kappa,\xi} \leq \int^{c_d}_{\ov\gamma} \4\eps^{-1/2} \4dJ_{\kappa,\xi}(\eps)
 + c_d^{-1/2} \int^{\xi}_{\kappa} \4g (t)\4 dt.\nn
\eeqa

Using partial integration we have by \eqref{eq:gammacd} and the
definition of $\bar{\gamma}$, \beqa\label{eq:I-estimate}
I_{\kappa,\xi} &\le& c_d^{-1/2}\underbrace{J_{\kappa,\xi}(c_d)}_{=
\;G(\kappa,\xi)} -\,
\bar{\gamma}^{-1/2}\underbrace{J_{\kappa,\xi}(\bar{\gamma})}_{= \;0}
+\ffrac{1}{2} \4 \int^{c_d}_{\ov{\gamma}}\4 \eps^{-3/2}
J_{a,b}(\eps) \, d\eps + \4 c_d^{-1/2}\4 G(\kappa,\xi)\nonumber\\
&=&\ffrac{1}{2} \4 \int^{c_d}_{\ov{\gamma}}\4 \eps^{-3/2}
J_{a,b}(\eps) \, d\eps + \4 2c_d^{-1/2}\4 G(\kappa,\xi) .\eeqa
Furthermore, \,$ {M}(t) \geq ( M_{1,t})^d \geq r^{-d}$ \, (see
\eqref{eq:minima0})  implies together with Lemma \ref{Mzerlegung}
 \beqa
J_{\kappa,\xi}(\eps) &\leq& \int_{\kappa}^{\xi}g(t)I_{\{\left(M_{1,t}\right)^d \leq \eps\}}\, dt =  \int_{\kappa}^{\xi}g(t)I_{\{M_{1,t} \leq \eps^{1/d}\}}\, dt \nonumber \\
&\leq& \int_{\kappa}^{\xi}g(t)I_{\{M^+_{1,t} \leq \eps^{1/d}\}}\,dt +\int_{\kappa}^{\xi}g(t)I_{\{ M^-_{1,t} \leq \eps^{1/d}\}}\, dt \nonumber \\[2mm] &=& H_{\kappa,\xi,\Q^+}(\eps^{1/d}) + H_{\kappa,\xi,\Q^-}(\eps^{1/d}),\label{eq:J-estimate}
\eeqa
where $H_{\kappa,\xi,\Q^{\pm}}$ is defined as in \eqref{eq:Hdef} in Lemma \ref{integral}. The smallest and the largest eigenvalue of $\Q^{\pm}$ is denoted by $q_0^{\pm}$ and  $q^{\pm}$, respectively.\\
Since $r^{-d}\leq \eps \leq c_d $ and hence $ r^{-1}\leq \eps^{1/d}\leq (2d)^{-1}\leq (2d^{\pm})^{-1}$ Lemma \ref{integral} can be applied and  by changing the variable $v=r^d\4 \eps $ we obtain\\

$
\int_{\bar{\gamma}}^{c_d} \eps^{-3/2}H_{\kappa,\xi,\Q^{\pm}}\left(\eps^{1/d}\right)d\eps  $
\beqa
\hskip1mm
&\leq&\! c(d^{\pm})
 \int_{\bar{\gamma}}^{c_d} \eps^{-3/2}\left(\ffrac{q^\pm}{q^\pm_0} \4\eps^{2/d}
G(\kappa(\eps^{1/d}\4  r),\xi)\, +\ffrac{1}{q_0^{\pm}}\ffrac{\eps^{1/d}}{r} g(\kappa(\eps^{1/d}\4 r))\right)d\eps\nn\\
&\leq& \!c(d^{\pm})\int_{\gamma}^{D}r^{\frac{3d}{2}-2} v^{-\frac{3}{2}+\frac{1}{d}}\left(\ffrac{q^\pm}{q^\pm_0} v^{1/d}
G(\kappa(v^{1/d} ),\xi)\, +\ffrac{1}{q_0^{\pm}}\, g(\kappa(v^{1/d}))\right)r^{-d}dv\nn\\
&=&\! r^{\frac{d}{2}-2}\cdot c(d^{\pm})\int_{\gamma}^{D} v^{-\frac{1}{2}+\frac{1}{d}}\left(\ffrac{q^\pm}{q^\pm_0} \4v^{1/d}
G(\kappa(v^{1/d} ),\xi) \, +\,  \frac{1}{q_0^{\pm}}g(\kappa(v^{1/d})) \right)\ffrac{dv}{v}.\nn
\eeqa
Analyzing the proof of Lemma \ref{integral} we may assume w.l.o.g that the constant $c(d)$ is monotone increasing in $d$. Since $ q_0 = \min\left\{q_0^+;q_0^-\right\} $ and $ q = \max\left\{q^+;q^-\right\} $, we have\\

$ 
\int_{\bar{\gamma}}^{c_d} \eps^{-3/2}H_{\kappa,\xi,\Q^{\pm}}\left(\eps^{1/d}\right)d\eps  $
\beqa
&\leq&
r^{d/2-2}\cdot c(d)\int_{\gamma}^{D} v^{-\frac{1}{2}+\frac{1}{d}}\left(\ffrac{q}{q_0} \4v^{1/d}
G(\kappa(v^{1/d} ),\xi) \, +\,  \frac{1}{q_0}g(\kappa(v^{1/d})) \right)\ffrac{dv}{v}.\nn\\[-2mm]&&\label{eq:H-estimate}
\eeqa
Thus we conclude by using \eqref{eq:I-estimate}, \eqref{eq:J-estimate},
and \eqref{eq:H-estimate}
\beqa
I_{\kappa ,\xi } &\ll_d& r^{d/2-2}\int_{\gamma}^{D} v^{-1/2+1/d}\left(\ffrac{q}{q_0} \4v^{1/d}
G(\kappa(v^{1/d} ),\xi) \, +\,  \frac{1}{q_0}g(\kappa(v^{1/d})) \right)\ffrac{dv}{v}\nn\\&& +\quad G(\kappa,\xi),\nn
\eeqa
which proves \eqref{eq:finalint}.
This completes the proof of Lemma \ref{final}.
\bwend

\begin{lemma}\label{irrational}
Let $0< \kappa<\xi <\infty$. Then
\beqa
 \lim\limits_{r \rightarrow \infty} \;\;\inf\limits_{t \in [\4\kappa,\xi\4]}
\bigl(r\4 M_{1,t}\bigr) \cdots \bigl(r\4 M_{d,t}\bigr) = \infty
\nonumber
\eeqa
provided that $\Q$ is irrational.
\end{lemma}

\bw \cite{goetze:2004}, p. 224, Lemma 3.11 or \cite{elsner:2006}, p. 47, Lemma 2.4.24 \bwend \vskip-10mm \enlargethispage{10mm}\phantom{...}
\bibliographystyle{alpha}\bibliography{ref}

\begin{thebibliography}{EMM98}

\bibitem[BG97]{bentkus-goetze:1997}
V.~Bentkus and F.~G{\"o}tze.
\newblock On the lattice point problem for ellipsoids.
\newblock {\em Acta Arith.}, 80(2):101--125, 1997.

\bibitem[BG99]{bentkus-goetze:1999}
V.~Bentkus and F.~G{\"o}tze.
\newblock Lattice point problems and distribution of values of quadratic forms.
\newblock {\em Ann. of Math. (2)}, 150(3):977--1027, 1999.

\bibitem[Cas78]{cassels:1978}
J.~W.~S. Cassels.
\newblock {\em Rational quadratic forms}, volume~13 of {\em London Mathematical
  Society Monographs}.
\newblock Academic Press Inc. [Harcourt Brace Jovanovich Publishers], London,
  1978.

\bibitem[Dav58]{davenport:1958}
H.~Davenport.
\newblock Indefinite quadratic forms in many variables. {II}.
\newblock {\em Proc. London Math. Soc. (3)}, 8:109--126, 1958.

\bibitem[DH46]{davenport-heilbronn:1946}
H.~Davenport and H.~Heilbronn.
\newblock On indefinite quadratic forms in five variables.
\newblock {\em J. London Math. Soc.}, 21:185--193, 1946.

\bibitem[DL72]{davenport-lewis:1972}
H.~Davenport and D.~J. Lewis.
\newblock Gaps between values of positive definite quadratic forms.
\newblock {\em Acta Arith.}, 22:87--105, 1972.

\bibitem[DM93]{dani-margulis:1993}
S.~G. Dani and G.~A. Margulis.
\newblock Limit distributions of orbits of unipotent flows and values of
  quadratic forms.
\newblock In {\em I. M. Gel$'$fand Seminar}, volume~16 of {\em Adv. Soviet
  Math.}, pages 91--137. Amer. Math. Soc., Providence, RI, 1993.

\bibitem[Els06]{elsner:2006}
G.~Elsner.
\newblock Distributions of values of indefinite forms and higher-order spectral
  estimates for finite markov chains.
\newblock {\em PhD thesis}, 2006.

\bibitem[EMM98]{eskin-margulis-mozes:1998}
A.~Eskin, G.~Margulis, and S.~Mozes.
\newblock Upper bounds and asymptotics in a quantitative version of the
  {O}ppenheim conjecture.
\newblock {\em Ann. of Math. (2)}, 147(1):93--141, 1998.

\bibitem[G{\"o}t04]{goetze:2004}
F.~G{\"o}tze.
\newblock Lattice point problems and values of quadratic forms.
\newblock {\em Invent. Math.}, 157(1):195--226, 2004.

\bibitem[Jar28]{jarnik:1928}
V.~Jarn\'{{\i}}k.
\newblock \"{U}ber {G}itterpunkte in mehrdimensionalen {E}llipsoiden.
\newblock {\em Math. Ann.}, 100(1):699--721, 1928.

\bibitem[Mar89]{margulis:1987}
G.~A. Margulis.
\newblock Discrete subgroups and ergodic theory.
\newblock In {\em Number theory, trace formulas and discrete groups (Oslo,
  1987)}, pages 377--398. Academic Press, Boston, MA, 1989.

\bibitem[Mar97]{margulis:1997}
G.~A. Margulis.
\newblock Oppenheim conjecture.
\newblock In {\em Fields Medallists' lectures}, volume~5 of {\em World Sci.
  Ser. 20th Century Math.}, pages 272--327. World Sci. Publishing, River Edge,
  NJ, 1997.

\bibitem[Mum83]{mumford:1983}
D.~Mumford.
\newblock {\em Tata lectures on theta. {I}}, volume~28 of {\em Progress in
  Mathematics}.
\newblock Birkh\"auser Boston Inc., Boston, MA, 1983.

\bibitem[Opp29]{oppenheim:1929}
A.~Oppenheim.
\newblock The minima of indefinite quaternary quadratic forms.
\newblock {\em Proc. Nat. Acad. Sci. USA}, 15:724--727, 1929.

\bibitem[Opp31]{oppenheim:1931}
A.~Oppenheim.
\newblock The minima of indefinite quaternary quadratic forms.
\newblock {\em Ann. of Math. (2)}, 32(2):271--298, 1931.

\end{thebibliography}
\end{document}